\documentclass[12pt,a4paper]{amsart}
\usepackage{amsmath}
\usepackage{amscd}
\usepackage{amssymb}

\numberwithin{equation}{section}
\theoremstyle{plain} %% This is the default
\newtheorem{thm}{Theorem}[section]
\newtheorem{cor}[thm]{Corollary}
\newtheorem{lem}[thm]{Lemma}
\newtheorem{prop}[thm]{Proposition}

\newtheorem{defn}[thm]{Definition}

\theoremstyle{remark}

\newcommand{\thmref}[1]{Theorem~\ref{#1}}
\newcommand{\propref}[1]{Proposition~\ref{#1}}

\newcommand{\lemref}[1]{Lemma~\ref{#1}}

\newsavebox{\SmallMathBox}

\def\pdo{\psi{\rm do}}

\def\Ci{C^\infty}

\def\dd{\partial}

\def\Di{D\kern -.65em /}
\def\Dii{D\kern -.45em /}
\def\di{{\dd}\kern -.55em /}
\def\dii{{\dd}\kern -.40em /}

\def\noi{\noindent}

\def\to{\rightarrow}
\def\too{\longrightarrow}

\def\mtoo{\longmapsto}

\def\re{{\rm Re}}

\def\pp{{\bf p}}

\def\qq{{\bf q}}

\def\Aa{{\mathcal A}}

\def\Cc{{\mathcal C}}
\def\Dd{{\mathcal D}}
\def\Ee{{\mathcal E}}

\def\Ii{{\mathcal I}}

\def\Ll{{\mathcal L}}

\def\Nn{{\mathcal N}}

\def\Rr{{\mathcal R}}

\def\={\cong}
\def\>{\supset}
\def\<{\subset}
\def\ii{^{-1}}
\def\si{^{-s}}

\def\pp{^{\perp}}
\def\12{\frac{1}{2}}
\def\2{\Dd}
\def\3{\Nn}
\def\4{\Rr}
\def\6{\cup}
\def\8{\otimes}
\def\0{^{\circ}}
\def\){\hfill{\ \qed}\enddemo}

\def\a{\alpha}

\def\b{\beta}
\def\C{\mathbb{C}}

\def\d{\delta}
\def\D{\Delta}
\def\e{\varepsilon}

\def\g{\gamma}
\def\G{\Gamma}

\def\th{\theta}

\def\la{\lambda}

\def\o{\infty}
\def\p{\pi}

\def\s{\sigma}
\def\Si{\Sigma}
\def\z{\zeta}

\def\ch{\mbox{\rm ch}}

\def\sdet{\mbox{\rm sdet\,}}

\def\End{\mbox{\rm End}}

\def\Ind{\mbox{\rm Ind\,}}

\def\LIM{\mbox{\rm LIM}}

\def\Ker{\mbox{\rm Ker}}

\def\ind{\mbox{\rm ind\,}}

\def\Si{S\kern -.65em /}

\def\str{\mbox{\rm Str}}

\def\Tr{\mbox{\rm Tr\,}}

\def\Af{\mathbb{A}}

\def\Cf{\mathbb{C}}

\def\Nf{\mathbb{N}}

\def\Ps{\mathbb{P}}

\def\Rf{\mathbb{R}}

\def\Zf{\mathbb{Z}}

\def\xx{\textbf{x}}

\def\ord{\mbox{\rm ord}}

\def\cc{\textsf{c}}

\def\As{\textsf{A}}
\def\Bs{\textsf{B}}

\def\Ds{\textsf{D}}

\def\Fs{\textsf{F}}

\def\Is{\textsf{I}}

\def\Ks{\textsf{K}}

\def\Ps{\textsf{P}}
\def\Qs{\textsf{Q}}
\def\Rs{\textsf{R}}

\def\Ts{\textsf{T}}

\def\Ws{\textsf{W}}

\begin{document}

\title[]{
Zeta Forms and the Local Family Index Theorem\vskip 5mm {\Small
{\rm Simon Scott}}}

\maketitle

%\tableofcontents

%*****************************************************************

\section{Introduction and Preliminaries}
Let $X$ be a $\Ci$ $n$-dimensional compact Riemannian manifold
without boundary and let $E = E^+ \oplus E^-$ be  a
$\mathbb{Z}_2$-graded (super) complex vector bundle over $X$. We
write $\G(X,E)$ for the space of $\Ci$ sections of $E$ and $\tau$
for the involution defining the induced $\mathbb{Z}_2$-grading.
The super (or $\mathbb{Z}_2$-graded) trace of a trace class
operator $a$ on $\G(X,E)$ is defined by $\str(a) = \Tr(\tau a)$.
Let $A$ and $F$ be classical (one-step polyhomogeneous)
pseudodifferential operators ($\pdo$s) acting on $\G(X,E)$.
Suppose that $A$ is of order $\nu\in\Rf$ and that $F$ is elliptic
of positive integer order $k$ and such that there is an angle
$\th$ for which the principal symbol $\s_k(F)(x,\xi)$ has no
eigenvalues on $R_{\th} = \{re^{i\th} \ | \ r > 0\}$. In this
situation, Grubb-Seeley \cite{GS1} show that as $\la\too \o$ in an
open sub-sector of $\Cf$ around $R_{\th}$ there is an asymptotic
expansion of the resolvent supertrace for $m > (n+\nu)/k$
\begin{equation}\label{e:resasymp1}
\str(A(F - \la I)^{-m}) \sim \sum_{j=0}^{\o} \a_j
(-\la)^{\frac{\nu + n - j}{k} - m} + \sum_{k=0}^{\o}
(\a_k^{'}\log(-\la) + \a_k^{''})(-\la)^{-k-m} \ .
\end{equation}
On the other hand, for $\re(s) > (n+\nu)/k$ the complex powers $A
F_{\th}\si$ are trace class and a generalized (super) zeta
function can be defined by
\begin{equation*}\label{e:superzeta}
\zeta_{\th}(A,F,s) = \str(AF_{\th}\si) \ .
\end{equation*}
When $A$ is the identity we write $\zeta_{\th}(F,s) :=
\zeta_{\th}(I,F,s)$. It is well known \cite{GS1,GS2,S} that the
expansion \eqref{e:resasymp1} is essentially equivalent to the
meromorphic extension $\zeta(A,F,s)|^{\rm mer}$ of $\zeta(A,F,s)$
(omitting the $\th$ subscript) to all of $\Cf$ with the
singularity structure
\begin{equation}\label{e:zetasing}
\G(s)\,\zeta(A,F,s)|^{\rm mer}  \sim \sum_{j=0}^{\o} \frac{a_j}{s
+ \frac{j + \nu - n}{k}} - \frac{\str(A\Pi_0(F))}{s} +
\sum_{k=0}^{\o} \frac{a_k^{'}}{(s + k)^2} + \frac{a_k^{''}}{s+k} \
,
\end{equation}
where $\Pi_0(F)$ is the orthogonal projection onto the kernel
$\Ker(F)$ of $F$. The coefficients in \eqref{e:zetasing} differ
from those in \eqref{e:resasymp1} by universal multiplicative
constants. The $a_j, a_j^{'}$ are local, being determined by
finitely many homogeneous terms of the local symbol expansions,
while the $a_j^{''}$ depend globally on $A, F$ and the bundle $E$.

Since $\G(s)\ii = s + o(s)$ around $s=0$, \eqref{e:zetasing}
implies that $\zeta_{\th}(A,F,s)|^{\rm mer}$ is holomorphic at
$s=0$ provided $a_0^{'} = 0$. In particular, this holds for the
zeta function $\zeta_{\th}(P^2,s)|^{\rm mer}$ associated to the
odd parity operator $P =
\begin{bmatrix}
  0 & P^- \\
  P^+ & 0 \\
\end{bmatrix}$, where $P^+ : \G(X,E^+)\too \G(X,E^-)$ is a classical
elliptic $\pdo$ of positive order, and $P^-$ its formal adjoint.
Since $P^+ P^-$ and $P^- P^+$ have identical non-zero spectrum
while $P^{-2s}$ vanishes on $\Ker(P)$ for $\re(s) >0$, it follows
that $\zeta_{\th}(P^2,s)|^{\rm mer}=0$. Evaluating at $s=0$ gives
the Atiyah-Bott-Seeley zeta function formula for the index
\begin{equation}\label{e:abs}
    \zeta_{\pi}(P^2, 0)|^{{\rm mer}} = 0 \ ;
\end{equation}
for, $\str(\Pi_0(P^2)) =  \dim\Ker(P^+) - \dim\Ker(P^-) :=
\ind(P)$ and hence from \eqref{e:zetasing} (and
\eqref{e:resasymp1}) equation \eqref{e:abs} is the identity
\begin{equation}\label{e:indexterms}
    \ind(P) = a_n + a_0^{''} = \a_n + \a_0^{''} \ .
\end{equation}
When $P$ is a differential operator, then $a_0^{''} = 0$ and
\eqref{e:indexterms} gives a formula for the index as the integral
over $X$ of a locally determined density.

Since $P^2$ is positive, \eqref{e:abs} and \eqref{e:indexterms}
are further equivalent for $t>0$ to the heat trace formula
$\ind(P) =  \str(e^{-tP^2})$; --- if $F$ is positive
\eqref{e:resasymp1} and \eqref{e:zetasing} are  equivalent
\cite{GS2} to a heat trace expansion  as $t \to 0+$ (with the same
coefficients as \eqref{e:zetasing})
\begin{equation}\label{e:heatracesasymp}
\str(Ae^{-tF}) \sim \sum_{j=0}^{\o} a_j t^{\frac{j-\nu-n}{k}} +
\sum_{k=0}^{\o} (-a_k^{'}\log t + a_k^{''})t^k \ .
\end{equation}

If $a_0^{'} = 0$, the next term up in the Laurent expansion of
$\zeta_{\th}(F,s)|^{{\rm mer}}$ around $s=0$ of a classical $\pdo$
$F$ is the logarithm of the regularized (graded- or super-)
zeta-determinant $\det_{\z,\th}F$. Thus
\begin{equation}\label{e:sdet}
\log{\rm det}_{\z,\th} F  = -
\frac{d}{ds}\zeta_{\th}(F,s)|_{s=0}^{{\rm mer}} =
\zeta_{\theta}(\log F , F, s)|_{s=0}^{{\rm mer}} \ .
\end{equation}
It is consistent to also write $\sdet_{\z,\th}F$ for the super
zeta-determinant, but here we prefer to retain the usual notation
unless we need to emphasize the grading. Notice that in the case
of the trivial grading the supertrace reduces to the usual
operator trace and so $\det_{\z,\th}F$ then coincides with the
usual ungraded zeta determinant, while, for example, for
even-parity $F = F^+ \oplus F^-$ one has
\begin{equation*}%\label{e:sdet}
{\rm sdet}_{\z,\th} F  = \frac{\det_{\z,\th} F^+ }{\det_{\z,\th}
F^- } \
\end{equation*}
with $\det_{\z,\th} F^{\pm}$ ungraded zeta determinants.

\vskip 1mm

In this paper we extend these constructions to geometric families
of $\pdo$s. We consider a $\Ci$ fibration $\pi : M \too B$ of
finite-dimensional manifolds with closed boundaryless fibre $M_z =
\pi\ii(z)$ equipped with a  Riemannian metric $g_{M/B}$  on the
vertical tangent bundle $T(M/B)$. Let $|\wedge_{\pi}| =
|\wedge(T^*(M/B))|$ be the line bundle of vertical densities,
restricting on each fibre to the usual bundle of densities
$|\wedge_{M_z}|$ along $M_z$. Let $\Ee = \Ee^+\oplus \Ee^-$ be a
vertical Hermitian super bundle over $M$ and let $\pi_*(\Ee) =
\pi_*(\Ee^+)\oplus\pi_*(\Ee^-)$ the graded infinite-dimensional
Frechet bundle with fibre
$\G(M_z,\Ee^z\otimes|\wedge_{M_z}|^{1/2})$ at $z\in B$, where
$\Ee^z$ is the super bundle over $M_z$ obtained by restriction of
$\Ee$. By definition, a $\Ci$ section of $\pi_*(\Ee)$ over $B$ is
a $\Ci$ section of $\Ee\otimes|\wedge_{\pi}|^{1/2}$ over $M$, and
more generally the de-Rham complex of $\Ci$ forms on $B$ with
values in $\pi_*(\Ee)$ is defined by
\begin{equation*}
\Aa(B,\pi_*(\Ee))  = \G(M,\pi^*(\wedge T^* B)\otimes
\Ee\otimes|\wedge_{\pi}|^{1/2})
\end{equation*}
with $\otimes$ the graded tensor product. We write $\Psi(\Ee)$ for
the infinite-dimensional bundle of algebras with fibre
$\Psi(\Ee^z)$ the space of classical $\pdo$s on
$\G(M_z,\Ee^z\otimes|\wedge_{M_z}|^{1/2})$. A section
$\Fs\in\Aa(B,\Psi(\Ee))$ defines a smooth family of $\pdo$s with
differential form coefficients parameterized by $B$. If the smooth
family of $\pdo$s
\begin{equation}\label{e:degreezero}
\Ps := \Fs_{[0]}\in \G(B,\Psi(\Ee))
\end{equation} defined by the form degree
zero component of $\Fs$ has positive order and admits a spectral
cut $R_{\th}$ then, for an auxiliary family of $\pdo$s
$\As\in\Aa(B,\Psi(\Ee))$, we use the fibrewise supertrace to
define for $\re(s)>>0$ a mixed degree differential form
\begin{equation}\label{e:superzetaform0}
\zeta_{\th}(\As,\Fs,s) = \str(\As \Fs_{\th}\si) \in \Aa(B) .
\end{equation}
When $\As = \Is$ is the vertical family of identity operators we
write $\zeta_{\th}(\Fs,s) := \zeta_{\th}(\Is,\Fs,s)$. In a similar
way to the single operator case, an analysis of the asymptotic
behavior of the corresponding resolvent trace differential form
defined for sufficiently large $m$ and $|\la|$
$$\str(\As(\Fs-\la\Is)^{-m})\in\Aa(B,\Psi(\Ee)) \ ,$$
shows that if the kernels $\Ker P_z$ of the family
\eqref{e:degreezero} have constant dimension, then the zeta trace
form \eqref{e:superzetaform0} extends meromorphically on $\Cf$ to
a form $\zeta_{\th}(\As,\Fs,s)|^{\rm mer}$.

For a family of strictly positive operators $\Fs$ one has
additionally the heat trace form
\begin{equation}\label{e:heatraceform}
\str(\As e^{-\Fs})\in \Aa(B) \ ,
\end{equation}
which, in the case when $\As = \Is$ and $\Fs$ is the curvature
form of a superconnection, is the object of interest in Bismut's
heat equation proof of the local Atiyah-Singer index theorem for
families of Dirac operators. Resolvent trace forms, or
equivalently the zeta trace forms, are, however, more general
geometric invariants which are also defined for families of
non-positive operators --- for example, the zeta form for a family
of self-adjoint first-order elliptic differential operators over
an odd-dimensional manifold. This is concordant with
\cite{GS1,GS2} where the resolvent trace and power operators were
shown to provide a considerably more powerful tool for $\pdo$
analysis, than heat kernel methods alone -- a principle which
applies equally well to families of $\pdo$s. These constructions
are for general elliptic families, the particular case of a family
of Dirac operators is considered when explicit formulas are sought
for the locally determined coefficients in the trace expansions,
such as in the case of the local index density and for such
geometric results heat kernel methods have so far been the most
effective. \vskip 3mm

It is useful, then, to use these methods to compute higher
geometric invariants as generalized $\zeta$-forms, such as
Wodzicki residue trace forms which extend the usual residue trace
functional in so far as it vanishes on super commutators of
families of $\pdo$s, the Kontsevich-Vishik trace form, eta-forms,
analytic torsion forms, and the extension of the corresponding
zeta form invariants to families of singular manifolds. Here we
illustrate the methods with an alternative view point onto the
Atiyah-Singer family index theorem.

\vskip 3mm

The following formulas generalize the Atiyah-Bott-Seeley formula
\eqref{e:abs}.
\begin{thm}\label{t:1}
Let $\Af^2\in\Aa(B,\Psi(\Ee))$ be the curvature form of a
superconnection $\Af$ on $\pi_*(\Ee)$ adapted to a smooth family
of elliptic $\pdo$s $\Ps = (P_z \ | \ z\in B)\in \G(B,\Psi(\Ee))$.
Suppose that the family of kernels $\Ker P_z = \Ker P^+_z \oplus
\Ker P^-_z$ have constant dimension as $z$ varies in $B$, forming
a super bundle $\Ker \Ps$ over $B$. Then $\z_{\pi}(\Af^2,s)|^{{\rm
mer}}$ is canonically exact in $\Aa(B)$. There is a canonical
transgression form $\tau_{\Af^2} \in \Aa(B)$ such that the
following formula holds in $\Aa(B)$
\begin{equation}\label{e:familyzetaindex}
\sum_{k=0}^{\dim B}\frac{1}{k!}\ \z_{\pi}(\Af^2,-k)|^{{\rm mer}} =
d\tau_{\Af^2}
\end{equation}
and implies the Family Index Theorem transgression formula for the
Chern character form. Precisely, replacing $\Af$ by the
$t$-rescaled superconnection $\Af_t$, the regularized limit as
$t\to 0+$ of \eqref{e:familyzetaindex} is the formula
\begin{eqnarray}
\ch(\Ker \Ps, \nabla^0) & =  & \LIM_{t \to 0} \left( \ch(\Af_t) -
d\int_{t}^{\o}\str(\dot{\Af}_{\s}e^{-\Af_{\s}^2})\, d\s
 \right) \nonumber\\
 & = & \int_{M/B}
\str(A(x)) - d\,\LIM_{t \to 0}\left(
\int_{t}^{\o}\str(\dot{\Af}_{\s}e^{-\Af_{\s}^2})\, d\s\right) \
,\label{e:familycherncharacterindex2}
\end{eqnarray}
where $A$ is a section of the bundle $\wedge
T^{*}M\otimes\End(\Ee)$ over $M$ whose d-form component $A_{[d]}$
is the coefficient of $\la^{-m-1-d/2}$ in the asymptotic expansion
of the kernel of the resolvent  $(\Af^2
-\la\Is)^{-m}\in\Aa(B,\Psi(\Ee))$ as $\la\to \o$ in a subsector of
$\Cf$.
\end{thm}

In \eqref{e:familycherncharacterindex2},  $\Af_t$ is defined by
multiplying the form degree $i$ component $\Af_{[i]}$ of $\Af$ by
$t^{1-i/2}$. The scaled super Chern character form is defined by
$\ch(\Af_t) = \str(e^{-\Af_{t}^2})$, while for a finite rank super
bundle $V$ with connection $\nabla$, $\ch(V,\nabla) =
\str(e^{-\nabla^2})$ is the classical graded Chern character form.
If $\Pi_0$ is the smooth family of smoothing operator projections
onto $\Ker \Ps$, $\nabla^0 = \Pi_0\cdot \Af_{[1]}\cdot\Pi_0$ is
the induced classical connection on $\Ker \Ps$. $\int_{M/B}:
\Aa^k(M)\too\Aa^{k-n}(B)$ denotes integration over the fibre,
while  the regularized limit $\LIM_{t\to 0}$ picks out the $t^0$
term in the asymptotic expansion as $t\to 0+$.

\vskip 5mm

Notice from \eqref{e:familyzetaindex} that it is now not just the
value of the zeta function at zero that determines the index, but
also its meromorphically continued value at a finite number of
negative integers.

\vskip 3mm

By standard index theory arguments \thmref{t:1} implies the
following general cohomological formula for any smooth family of
elliptic $\pdo$s.

\begin{cor} For any elliptic family of $\pdo$s $\Ps\in
\G(B,\Psi^{>0}(\Ee))$
 there exists a smooth family of smoothing operators
$\Ks\in \G(B,\Psi^{-\o}(\Ee))$ such that $\Ps + \Ks$ has constant
kernel dimension and $\Ind(\Ps + \Ks) = \Ind(\Ps)$ in $K(B)$. The
following formula holds in $H^*(B)$
\begin{equation}\label{e:zetafamilycohomologyindex}
\sum_{k=0}^{\dim B}\frac{1}{k!}\ \z_{\pi}(\Af^2 + \Ks,-k)|^{{\rm
mer}} = 0
\end{equation}
and implies the cohomological family index theorem in $H^*(B)$
\begin{equation}\label{e:cohomologyfamilycherncharacterindex}
\ch(\Ind \Ps) = \LIM_{t \to 0}\ch(\Af_t) = \int_{M/B} \str(A(x))
\ .
\end{equation}
\end{cor}

Equivalently, this zeta form  may be viewed as a map $K(B) \too
H^*(B)$ which vanishes identically -- naturally generalizing
\eqref{e:abs}. More generally, the meromorphically continued zeta
form for a smooth family of $\pdo$s is a K-theoretic invariant.

\vskip 2mm

Following the single operator case, for a general smooth family of
$\z$-admissible  $\pdo$s $\Fs\in\Aa(B,\Psi(E))$ with spectral cut
$R_{\th}$ the next term up in the Laurent expansion around $s=0$
of $\zeta_{\th}(\Fs,s)$ defines the logarithm of the super zeta
determinant form
\begin{equation}\label{e:zdetform}
    {\rm det}_{\z,\th}\Fs \in \Aa(B) \ .
\end{equation}
This is a non-local mixed-degree differential form invariant which
extends the classical operator zeta determinant. We may on
occasion write $\sdet_{\z,\th}\Fs$ for \eqref{e:zdetform} to
emphasize the grading.

\begin{lem}\label{p:zetaperturbation}
Let $\Ps = \Fs_{[0]}\in \G(B,\Psi^{>0}(\Ee))$ be the degree zero
component of $\Fs$. Then $\Ps$ is $\z$-admissible and
\begin{equation}\label{e:zetaperturbation}
    {\rm det}_{\z,\th}\Fs  = {\rm det}_{\z,\th}\Ps +
    \omega_{\z,\th}(\Fs) \ ,
\end{equation}
where $ \omega_{\z,\th}(\Fs) \in\Aa^{>0}(B) = \sum_{k\geq
1}\Aa^k(B)$.
\end{lem}
Hence the degree zero part of ${\rm det}_{\z,\th}\Fs $ coincides
with the classical zeta determinant function $z\mtoo \det_{\z,\th}
P_z$. \lemref{p:zetaperturbation} is proved in Section 2.

\vskip 2mm

We use these methods to give the following geometric application
of the $\zeta$-determinant form. The total Chern class on the
semi-group ${\rm Vect}(X)$ of finite rank complex vector bundles
defines a stable characteristic class and so descends to a ring
homomorphism
\begin{equation*}
    \cc : K(X) \too H^{*}(B) \ ,
\end{equation*}
on the K-ring of virtual bundles, which is an isomorphism on the
rational coefficient ring. For an element $[V^+] - [V^-] \in
K(B)$, represented by $V^{\pm}\in {\rm Vect}(X)$, one has
\begin{equation}\label{e:virtualchern}
    \cc([V^+] - [V^-]) = \frac{\cc(V^+)}{\cc(V^-)} \ .
\end{equation}
To construct a de-Rham representative for the cohomology identity
\eqref{e:virtualchern} we may use Quillen's observation \cite{Q}
that just as form representatives for the characteristic classes
of a vector bundle can be computed from a connection, so the
characteristic class forms of a virtual bundle can be computed
from a superconnection, extending Chern-Weil theory from ${\rm
Vect}(X)$ to $K(X)$. This applies to the Chern class in the
infinite-dimensional setting in the following way.

%For a superconnection $\Af$ on a finite-rank complex graded vector
%bundle $V = V^+\oplus V^-$, the total (super) Chern form is
%defined by the super determinant form
%\begin{equation*}
%\cc(V,\Af) = {\rm sdet}(I + \Af^2) :=\exp(\str(\log(I + \Af^2)))
%\in \Aa(X) \ .
%\end{equation*}
%By a straightforward argument (generalized in Section 3 to the
%infinite-dimensional case) $\cc(V,\Af)$ is a closed differential
%form and a homotopy invariant of the superconnection with $
%\cc([V^+] - [V^-])  = \cc(V,\Af)$ in $H^*(X)$. In particular,
%$\cc(V,\Af)$ is a homotopic form to $\cc(V^+,\nabla^+)/ \cc(V^-,
%\nabla^-)$, where $\nabla^{\pm}$ are connections on $V^{\pm}$ and
%$\cc(V^{\pm},\nabla^{\pm}) = \det(I + (\nabla^{\pm})^2) :=
%\exp(\tr\log(I + (\nabla^{\pm})^2))$ the classical (ungraded)
%Chern forms (we have omitted the usual normalizing factor of
%$i/2\pi$).

%In the infinite-dimensional situation this generalizes as follows.

\begin{thm}\label{t:2}
Let $\Af$ be a superconnection adapted to a family of self-adjoint
Dirac operators $\Ds = \begin{bmatrix}
  0 & D^- \\
  D^+ & 0 \\
\end{bmatrix}$ associated to a Clifford connection over a Riemannian
fibration of spin manifolds $\pi : M \to B$. The zeta-Chern form
defined by the super zeta determinant form
\begin{equation}\label{e:zetachern}
    \cc_{\, \z}(\Af) = {\rm sdet}_{\z,\pi}(\Is + \Af^2) \in \Aa(B)
\end{equation}
is a closed differential form and a homotopy invariant of $\Af$
representing the Chern class $\cc(\Ind(\Ds))\in H^*(B)$ of the
index bundle.

For $t>0$, let $\Af_t$ be the scaled superconnection. If $\Ds$ has
constant kernel dimension, defining a super bundle $\Ker(\Ds)\to
B$, then there is a Chern-Simons generalized zeta form
$\omega_{t,\o}(\Af)\in\Aa(B)$ such that for all $t>0$ the
transgression formula
\begin{equation}\label{e:cherntransgression}
    \cc(\Ker (\Ds),\nabla^0) = \cc_{\, \z}(\Af_t) +
    d\omega_{t,\o}(\Af)
\end{equation}
holds in $\Aa(B)$, where $\nabla^0$ is defined in \thmref{t:1}.
The left-side of \eqref{e:cherntransgression} is  the classical
(super) Chern form, defined for  a finite-rank complex graded
vector bundle $V$ with connection $\nabla$ by \, $\cc(V,\nabla) =
{\rm sdet}(I + \nabla^2) :=\exp(\str(\log(I + \nabla^2)))$. (A
precise formula for $\omega_{t,\o}(\Af)$ is given in Section 5.)

When $\Af_t$ is the Bismut superconnection then the form $\cc_{\,
\z}(\Af_t)$ has a limit in $\Aa(B)$ as $t\to 0+$ with
\begin{equation}\label{e:chernlimit}
\lim_{t\to 0+} \cc_{\, \z}(\Af_t)  = \prod_{j=1}^{[\dim
B/2]}e^{(-1)^{j-1} (j-1)! \left((2\pi i)^{-n}\int_{M/B}
\widehat{A}(M/B) \ \ch^{'}(\Ee) \right)_{[2j]}} \ ,
\end{equation}
where  $\widehat{A}(M/B)  =
\det^{1/2}\left(\frac{R^{M/B}/2}{\sinh(R^{M/B}/2)}\right)$ is the
vertical $\widehat{A}$-genus form and   $\ch^{'}(\Ee)$ is the
twisted Chern character form for $\Ee$. Here $\tau_{[j]}$ is the
$j$-form component of $\tau\in\Aa(B)$.
\end{thm}

By standard index theory arguments the constant kernel condition
can be dropped for the cohomological formula for $\cc(\Ind(\Ds))$:
\thmref{t:2} implies the following Chern class Families Index
Theorem.
\begin{cor}
For any smooth family of Dirac operators $\Ds$ associated to a
Clifford connection one has in $H^*(B)$
\begin{equation}\label{e:chernAtiyahSingerindex}
\cc(\Ind(\Ds))  = \prod_{j=1}^{[\dim B/2]}e^{(-1)^{j-1} (j-1)!
\left((2\pi)^{-n}\int_{M/B} \widehat{A}(M/B) \ \ch^{'}(\Ee)
\right)_{[2j]}} \ .
\end{equation}
\end{cor}

\noi {\bf Remarks.}  \newline \noi {\rm [1]}  As a generalized
zeta-determinant, the form $\cc_{\, \z}(\Af)$ is a highly
non-local invariant. That the identities \eqref{e:chernlimit},
\eqref{e:chernAtiyahSingerindex} hold is due to a localization of
$\cc_{\, \z}(\Af_t)$ when $\Af_t$ is the scaled Bismut
superconnection
\begin{equation*}
   \cc_{\, \z}(\Af_t) = \cc_{{\rm local},\z}(\Af) + \cc_{{\rm
   global},\z}(\Af)(t) + d\omega_t
\end{equation*}
into a local term plus a global term $\cc_{{\rm
global},\z}(\Af)(t)$ which is $O(t^{1/2})$  and an exact global
term which is $O(1)$ as $t\to 0+$. For a general superconnection
the form $\cc_{\, \z}(\Af_t)$ does not converge as $t\to 0+$.

\vskip 2mm

\noi {\rm [2]} The determinant line bundle $\Ll\to B$ of a family
of elliptic $\pdo$s $\Ps\in\G(B,\Psi(\Ee))$ of odd parity is the
complex line bundle with fibre $\wedge^{{\rm
max}}\Ker(P^+_z)^*\otimes\wedge^{{\rm max}}\Ker(P^-_z)$. The
Quillen-Bismut-Freed connection is defined on $\Ll$ via a
meromorphically continued zeta trace \cite{Q,BF} and one has:
\begin{cor}
The curvature of the Quillen-Bismut-Freed connection on $\Ll$ is
\begin{equation*}
R^{(\Ll)} = \LIM_{t\to 0+}\cc_{\, \z}(\Af_t)_{[2]} \ .
\end{equation*}
\end{cor}

\section{Asymptotic Expansion of the Resolvent Trace Form}

Let $\pi: M\too B$ be a smooth family of closed Riemannian
manifolds with Hermitian vector bundle $\Ee\too M$, as in the
introduction. Let $\Psi(\Ee)$ be the bundle of subalgebras of
$\End(\pi_*(\Ee))$ of classical $\pdo$ s, and let
$\Psi^{\nu}(\Ee)$ (resp. $\Psi^{<\nu}(\Ee)$) be the subbundle of
operators of order $\nu\in \G(B)$ (resp. less than $\nu$). Thus
the fibre of $\Psi^{\nu}(\Ee)$ at $z\in B$ is the algebra
$\Psi^{\nu(z)}(M_z,\Ee_z)$ of classical $\pdo$s on $\G(M_z,\Ee_z)$
of order $\nu(z)\in\Rf$, and its sections are families of $\pdo$s
which in any local trivialization of $M$ and $\Ee$ over an open
subset $U\subset B$ depend smoothly on the local coordinates.

%$\Qs$ defines a classical
%$\pdo$ on $\G(M,\pi^*(\wedge T^* B)\otimes\Ee)$ which
%supercommutes with the action of $\Aa(B)$.
The fibre product $M\times_{\pi}M$ is the fibration over $B$ with
fibre $M_z\times M_z$ and vertical bundle $\Ee\boxtimes\Ee :=
p_1^*(\Ee)\otimes p_2^*(\Ee)^{*}$, where $p_1,p_2 : M\times_{\pi}M
\to M$ are the canonical projection maps. For a smooth family of
$\pdo$s with differential form coefficients $\Qs\in
\Aa(B,\Psi(\Ee))$, if $x\in M$ is not in the support of
$\psi\in\Aa(B,\pi(\Ee))$ then there is a smooth family of smooth
kernels on $M\times_{\pi}M \backslash \Delta(M)$, where
$\Delta(M)$ is the diagonal $\{(x,x) \ | \ x\in M\}$ in
$M\times_{\pi}M$,
\begin{equation*}
K(\Qs)\in \G(M\times_{\pi}M \backslash \Delta(M), \pi^*(\wedge
T^*B)\otimes(\Ee\otimes |\Lambda_{\pi}|^{1/2})\otimes(\Ee^*\otimes
|\Lambda_{\pi}|^{1/2}))
\end{equation*}
with $|\Lambda_{\pi}|^{1/2}$ the line bundle of half-densities
along the fibres of $M$, such that
\begin{equation}\label{e:Qkernel}
    (\Qs \psi)(x) = \int_{M/B} K(\Qs)(x,y) \psi(y) \ .
\end{equation}
Restricted to the fibre $M_z$ \eqref{e:Qkernel} reduces  to the
usual pointwise kernel formula; for $x\in M_z$ not in the support
of $\psi\in \G(M_z,\Ee_z)$
\begin{equation}\label{e:Qkernelz}
    (Q^z \psi)(x) = \int_{M_z} K(Q^z)(x,y) \psi(y) \ .
\end{equation}

For a family $\Qs\in \Aa(B,\Psi^{<-n}(\Ee))$ of order less than
the fibre dimension $K(\Qs)$ extends continuously across
$\Delta(M)$, and hence applying to $\Qs$ the $\Aa(B)$ valued
supertrace
\begin{equation*}
\str: \Aa(B,\Psi^{<-n}(\Ee)) \too \Aa(B) \ ,
\end{equation*}
defined fibrewise for $z\in B$ by the operator supertrace on
$\Psi^{<-n}(\Ee^z)$, defines a differential form
$\str(\Qs)\in\Aa(B)$. On the other hand, the restriction of
$K(\Qs)$ to the diagonal $M\subset M\times_{\pi}M$  defines a
continuous section of $\pi^*(\wedge T^*B)\otimes\End(\Ee)\otimes
|\Lambda_{\pi}|$ depending smoothly on the base parameters, and so
the $\End(\Ee)$-supertrace  defines a section $\str(K(\Qs)(x,x))$
in $\G(M,\pi^*(\wedge T^*B)\otimes |\Lambda_{\pi}|)$ which can be
integrated over the fibres and we have
\begin{equation}\label{e:str}
    \str(\Qs) = \int_{M/B} \str\left(K(\Qs)(x,x)\right) \in \Aa(B) \ .
\end{equation}

For the general case, $\Qs\in \Aa(B,\Psi(\Ee))$ means that in any
local trivialization of $M$ and $\Ee$ the operator $\Qs$ is
represented by a {\em vertical polyhomogeneous symbol}. This means
the following. Assume a local trivialization $M_{|U_B} \cong
U_B\times M_{z_0}$ over an open subset $U_B\subset B$ with $z_0\in
U_B$, and a trivialization $\Ee_{U_M} \cong U_M \times \Rf^N$,
where $\Rf^N$ inherits a grading from $\Ee_{U_M}$, over an open
subset $U_M \subset \pi\ii(U_B)$. $U_M$ may be identified as a
product $U_M \cong U_B \times U_{z_0} \cong \Rf^{\dim B} \times
\Rf^n$ with $U_{z_0} = U_M \cap M_{z_0}$. A vertical symbol may be
written according to form degree as ${\bf q} = {\bf q}_{[0]} +
\ldots + {\bf q}_{[\dim B]}$ with
\begin{equation}\label{vertsymb}
{\bf q}_{[k]}(z,x,\xi)\in\G\left( \, U_B\times U_{z_0} \times
\Rf^{n}\backslash\{0\}, \ \pi^*(\wedge^k T^*
U_B)\otimes\Rf^N\times(\Rf^N)^*\right) \ ,
\end{equation}
where $\xi$ may be identified with an element of the vertical (or
fibre) cotangent space $T^{*}_x(M/B)$. Each $ {\bf q}_{[k]}$ can
be written (locally) as a finite sum of terms of the form
$\sum_{j=0}^{J}\omega_j\otimes  q_{[k],j}$ with
$\omega_j\in\Aa^k(U_B)$ a basis of local $k$-forms
 and $q_{[k],j}$ is a symbol (in the usual single manifold
 sense) of form degree zero. For clarity of exposition we shall
 assume for the moment that the vertical symbols
 are {\it simple}, meaning that they have the local form
 ${\bf q}_{[k]} = \omega_k\otimes q_{[k]}$ with just
 one term in each form degree ($J=0$). That ${\bf q}_{[k]}$ be a vertical (simple) symbol
 of order $\nu\in\G(B,\Rf^{\dim B + 1})$ is the growth requirement in fibre direction that
 for $k=0, \ldots, \dim B$ and for all multi-indices $\a,\b,\g$
\begin{equation}\label{e:verticalschwartzclass}
|\partial^{\g}_z\partial^{\a}_x\partial^{\b}_{\xi}{\bf
q}_{[k]}(z,x,\xi)| < C(1 + |\xi|)^{\nu_{k}(z) - |\b|} \ ,
\end{equation}
where $x\in U_{z_0}$, $z\in U_B$, and $\nu(z) =
(\nu_{1}(z),\ldots,\nu_{{\rm dim}(B)+1}(z))$, while on the
left-side $| \, . \, |$ is a choice of norm on $\pi^*(\wedge^k T^*
U_B)\otimes\Rf^N\times(\Rf^N)^*$. The estimate
\eqref{e:verticalschwartzclass} holds uniformly in $\xi$, and on
compact subsets of $U_B\times U_{z_0}$ uniformly in $(z,x)$. A
vertical symbol ${\bf q}$ is {\it classical} (1-step
polyhomogeneous) of order $\nu\in\G(B,\Rf^{\dim B + 1})$ if there
is an asymptotic expansion
\begin{equation}\label{e:qasymp}
{\bf q}_{[k]}(z,x,\xi) \sim \sum_{j\geq 0} {\bf
q}_{[k],j}(z,x,\xi)
\end{equation}
as $|\xi|\to \o$, meaning ${\bf q}_{[k]}(z,x,\xi) - \sum_{j
=0}^{N-1} {\bf q}_{[k],j}(z,x,\xi)$ is a symbol of order
${\nu_{k}(z) - N}$, and where ${\bf q}_{[k],j}(z,x,\xi)$ is
positively homogeneous of degree $\nu_k(z)-j$ in $\xi$, meaning
that for $t>0$
\begin{equation*}
{\bf q}_{[k],j}(z,x,t\xi) = t^{\nu_{k}(z)-j}{\bf
q}_{[k],j}(z,x,\xi) \ .
\end{equation*}
Then $\Qs: \Aa(B,\pi_{*}(\Ee))\to \Aa(B,\pi_{*}(\Ee))$  is a {\it
simple} vertical classical family of $\pdo$s parameterized by $B$
and of order $\nu\in\G(B,\Rf^{\dim B + 1})$, and we write $\Qs\in
\Aa(B,\pi_{*}(\Ee))$, if in any local trivialization, as above,
there is a simple vertical classical symbol ${\bf q} = {\bf
q}_{[0]} + \ldots + {\bf q}_{[\dim B]}$ of order $\nu$ such that
for $\psi$ with support in a compact subset of $U_M$
\begin{equation*}
\Qs\psi(x)  = \frac{1}{(2\pi)^n}\int_{\Rf^n}\int_{U_{z_0}}
e^{i(x-y).\xi} {\bf q}(z,x,\xi) \psi(y) \, dy\, d\xi  \ + \
\Rs\psi(x),
\end{equation*}
where $\Rs$ is a smooth family of smoothing operators, meaning
that it is defined by a smooth vertical kernel $K(\Rs)(x,y)$.

In particular, if as before $\Qs\in \Aa(B,\Psi^{<-n}(\Ee))$,
meaning that $\nu_i(z) < -n$ for $i=0,1,\ldots,\dim B+1$, then for
the corresponding vertical volume form, \eqref{e:str} is pointwise
the differential form in $\Aa(B)$
\begin{equation*}
\str(\Qs)(z)  = \int_{M/B}\frac{1}{(2\pi)^n}\int_{\Rf^n} {\rm
\str}({\bf q}(z,x,\xi))\, d\xi \, {\rm vol}_{M/B} \ .
\end{equation*}

Notice that writing ${\bf q}$ according to form degree corresponds
to writing
\begin{equation*}
\Qs = \Qs_{[0]} + \Qs_{[1]} + \ldots + \Qs_{[\dim B]}
\end{equation*}
where $\Qs_{[k]}\in \Aa^k(B,\Psi(\Ee))$ raises form degree in
$\Aa(B,\pi_*(\Ee))$ by $k$.  With respect to a local (weak)
trivialization of $\pi_*(\Ee)$ over $U\subset B$ one has
$$\Aa(U,\pi_*(\Ee)_{|U}) \cong \Aa(U)\otimes
\G(M_{z_0},\Ee^{z_0})$$ for $z_0\in U$ and so a general
$\Qs_{[k]}$ can be written locally as a sum of simple vertical
$\pdo$s $\Qs_{[k]|U} =  \sum_{j=0}^J \omega_{k,j} \otimes Q_{j} \
,$ with $\omega_{k,j}\in \Aa^k(U)$ and
$Q_{j}\in\G(U,\Psi^{\nu_{j}(z_0)}(\Ee_{z_0}))$.

Composition defines a canonical algebra structure on
$\Aa(B,\Psi(\Ee))$, coinciding with the usual pointwise structure
on $\Psi(\Ee_z)$, such that
\begin{equation}\label{e:multpdos}
\Aa^i(B,\Psi^{\nu}(\Ee)) \times \Aa^j(B,\Psi^{\mu}(\Ee)) \too
\Aa^{i+j}(B,\Psi^{\nu + \mu}(\Ee)) \ .
\end{equation}
The multiplication \eqref{e:multpdos} is defined locally  at the
symbol level;  if $\As\in \Aa^i(B,\Psi^{\nu}(\Ee)),  \Bs \in
\Aa^j(B,\Psi^{\nu}(\Ee))$ with simple local symbols given over
$U_M$ by
$${\bf a} = \omega_{[i]}\otimes a \in\G\left( \, (U_M\times_{\pi} U_M) \times
\Rf^{n}\backslash\{0\}, \ \pi^*(\wedge^i T^*
U_B)\otimes\Rf^N\times(\Rf^N)^*\right) \ ,$$ $${\bf b} =
\sigma_{[j]}\otimes b \in\G\left( \, (U_M\times_{\pi} U_M) \times
\Rf^{n}\backslash\{0\}, \ \pi^*(\wedge^j T^*
U_B)\otimes\Rf^N\times(\Rf^N)^*\right)$$ then $\As\Bs\in
\Aa^{i+j}(B,\Psi^{\nu + \mu}(\Ee))$ is defined by the vertical
polyhomogeneous symbol
\begin{equation}\label{e:aob}
{\bf a}\circ {\bf b} = \omega_{[i]}\wedge \sigma_{[j]}\otimes
(a\circ b) \ ,
\end{equation}
where, as elsewhere, $\otimes$ means the graded tensor product,
and $a\circ b \sim \sum_j (a\circ b)_j$ with
\begin{equation*}%\label{e:symprod}
(a\circ b)_j = \sum_{|\a| + k + l=j}\frac{(-i)^{\a}}{a
!}\partial^{\a}_{\xi}(a)_k
\partial^{\a}_{x}(b)_l \ .
\end{equation*}

The crucial property of a family of $\pdo$s in $\Aa(B,\Psi(\Ee))$
is that ellipticity properties are determined by its form degree
zero component. For a family of $\pdo$s $\Ps\in\G(B,\Psi(\Ee) =
\Aa^0(B,\Psi(\Ee))$, with differential form degree zero, the
principal symbol, defined in any local trivialization to be the
leading term ${\bf p}_{0}$ in the asymptotic expansion
\eqref{e:qasymp}, has an invariant realization as a smooth section
\begin{equation*}
{\bf p}_{0}\in \G(T(M/B),\varphi^*(\End(\Ee))) \ ,
\end{equation*}
where $\varphi : T(M/B) \too M$ is the tangent bundle along the
fibres.
\begin{defn}
A smooth family of $\pdo$s $\Qs\in \Aa(B,\Psi(\Ee))$ with
differential form coefficients is said to be elliptic with
principal angle $\th$ if its form degree zero component $\Ps :=
\Qs_{[0]}\in \G(B,\Psi(\Ee))$ is elliptic with principal angle
$\th$. This means that
\begin{equation*}
{\bf p}_{0} - \la \Ii\in \G(T(M/B)\backslash\{0\},p^*(\End(\Ee)))
\end{equation*}
is an invertible bundle map for $\la\in R_{\th} = \{re^{i\th} \ |
\ r > 0\}$, where $\Ii$ is the identity bundle operator,
\end{defn}

\begin{prop}\label{p:resolventformula}
If $\Qs\in \Aa(B,\Psi(\Ee))$ is elliptic with principal angle
$\th$, then there is an open sector $\G_{\th}\subset \Cf - \{0\}$
containing $R_{\th}$ such that on any compact codimension zero
submanifold $B_c$ of $B$ for large $\la\in \G_{\th}$ there is a
smooth family of vertical resolvent $\pdo$s
\begin{equation}\label{e:resolventexists}
(\Qs - \la \Is)\ii \in \Aa(B_c,\Psi(\Ee)) \ .
\end{equation}
Here $\Is$ is the vertical identity, defined by the symbol ${\bf
I} = 1\otimes I$, coinciding pointwise with the identity $I_z$ on
$\G(M_z,\Psi(\Ee_z))$, while the form degree zero component of
\eqref{e:resolventexists} coincides pointwise with the usual
$\pdo$ resolvent for $\Ps := \Qs_{[0]}$. Precisely, let
$$ \Qs_{[>0]} = \Qs - \Ps \in \Aa^{1}(B,\Psi(\Ee))\,)$$
be the component of $\Qs$ with non-zero form degree. Then for
large $\la\in \G_{\th}$ the following identity holds in
$\Aa(B_c,\Psi(\Ee))$
\begin{equation}\label{e:resolventformula}
(\Qs - \la \Is)\ii = (\Ps - \la \Is)\ii  + \sum_{k=1}^{\dim B}
(-1)^{k}(\Ps - \la \Is)\ii \left(\Qs_{[>0]}(\Ps - \la
\Is)\ii\right)^k \ .
\end{equation}
\end{prop}
\begin{proof}
The resolvent for $\Ps \in \Aa^0(B,\Psi^{\nu_0}(\Ee))$ can be
constructed using a standard procedure \cite{S,Sh}. Locally, with
respect to trivializations, for
$$\mu\in\G_{\th,m} =\{\upsilon\in \Cf\backslash\{0\} \ | \ \upsilon^m\in
\G_{\th}\}$$ inductively define vertical symbols ${\bf b}_j[\mu]
 \in\G\left( \, (U_M\times_{\pi} U_M) \times
\Rf^{n}\backslash\{0\}, \ \Rf^N\times(\Rf^N)^*\right)$ homogeneous
in $(\mu,\xi)$ of degree $-\nu -j$  by
$${\bf b}_0[\mu](z,x,\xi) = ({\bf
p}_{0}(z,x,\xi) - \mu^m I)\ii \ ,$$
\begin{equation}\label{e:parametrix}
{\bf b}_j[\mu](z,x,\xi) = ({\bf p}_{0}(z,x,\xi) - \mu^m I)\ii
\sum_{\stackrel{|\a| + k + l = j}{k<j}}\frac{(-i)^{\a}}{\a
!}\partial^{\a}_{\xi}{\bf p}_k (z,x,\xi)
\partial^{\a}_{x}{\bf b} _l [\mu](z,x,\xi) \ ,
\end{equation}
so that
\begin{equation}\label{e:parametrix2}
\left(\sum {\bf p}_k(z,x,\xi) - \mu^m I\right)\circ\left(\sum {\bf
b}_j[\mu]\right) \sim {\bf I} \ .
\end{equation}
It follows  that for $\la = \mu^m$ there exists a vertical
polyhomogeneous symbol ${\bf b}[\la] \sim \sum {\bf b}_j[\mu]$. By
a standard partition of unity construction, and using the local
trivializations, we can patch together to define $\Bs = {\rm
OP}({\bf b})\in\G(B,\Psi^{-\nu_0}(\Ee))$ and from
\eqref{e:parametrix2} $(\Ps-\la\Is)\Bs = I - \Rs$ with $\Rs\in
\G(B,\Psi^{-\o}(\Ee))$. The $L^2$ operator norm of $\Rs(z)$ is
$O(|\la|\ii)$ uniformly on $B_c$  and so $I-\Rs$ is invertible in
$\G(B_c,\Psi(\Ee))$ for sufficiently large $|\la|$ with
\begin{equation}\label{e:inverse}
(\Ps-\la\Is)\ii = \Bs + \Bs\sum_{j\geq 1}\Rs^j \ .
\end{equation}

The extension to \eqref{e:resolventformula} is a consequence of
the nilpotence of the differential form coefficients of
$\Qs_{[>0]} = \sum_{k=1}^{\dim B}\Qs_{[k]}$. From
\eqref{e:multpdos} and \eqref{e:aob}, if $$\As \in
\Aa^{>0}(B,\Psi(\Ee)) := \sum_{i=1}^{\dim B}\Aa^{i}(B,\Psi(\Ee)) \
,$$ so $\As$ (strictly) raises form degree in $\Aa(B,\pi_*(\Ee))$,
we have $\As^k = 0$ for $k > \dim B$. Hence the Neumann expansion
exists, consisting of only finitely many terms giving the identity
in $\Aa(B,\Psi(\Ee))$
\begin{equation*}
    (I - \As)\ii = I + \As + \ldots + \As^{\dim B} \ .
\end{equation*}
Since $(\Qs - \la \Is)\ii = (\Ps - \la \Is)\ii \left(\Is +
\Qs_{[>0]}(\Ps - \la \Is)\ii\right)$ and $\Qs_{[>0]}(\Ps - \la
\Is)\ii\in \Aa^{>0}(B,\Psi(\Ee))$ we reach the conclusion.
\end{proof}

\noi {\bf Remarks.}  \newline \noi {\rm [1]} Evidently, the proof
can also be carried out directly by constructing a parametrix for
the full local symbol ${\bf q}\in\G\left( \, (U_M\times_{\pi} U_M)
\times \Rf^{n}\backslash\{0\}, \ \pi^*(\wedge T^*
U_B)\otimes\Rf^N\times(\Rf^N)^*\right)$  from the vertical symbol
parametrix ${\bf b}$ for ${\bf p}(z,x,\xi) - \mu^m I$, as above,
and then using the differential form nilpotence of $({\bf q} -\la
I) - ({\bf p} -\la I)$ and the corresponding Neumann expansion,
defined for the symbol product \eqref{e:aob}, to construct a local
symbol parametrix for $({\bf q} -\la I)$ of the form
\eqref{e:resolventformula}.

\vskip 1mm

\noi {\rm [2]} In the following we shall for brevity not
distinguish between $B_c$ and $B$. In particular, this distinction
in \thmref{t:1} and \thmref{t:2} is not pertinent.

\vskip 3mm

\begin{thm}\label{t:0}
Let $\Fs=\sum_{k=0}^{\dim B}\Fs_{[k]}\in\Aa(B,\Psi(\Ee))$ be a
smooth simple family of elliptic $\pdo$s of constant order
$(v_0,v_2,\ldots,v_{\dim B})$, so that $\Ps = \Fs_{[0]}$ is an
elliptic family with parameter $\la\in \G_{\th}$ of constant order
$r=v_0>0$, and $v_j = \ord(\Fs_{[j]})\in\Rf^1$. Then with $w =
\sum_{k=1}^{\dim B+1}v_j$, if $\la\in \G_{\th}$ and $m > \frac{w +
n}{r}$ the resolvent derivative
\begin{equation}\label{e:resolventpowers}
\partial^{m-1}_{\la}(\Fs - \la \Is)\ii \in \Aa(B,\Psi^{< -n}(\Ee))
\end{equation}
is a smooth family of $\pdo$s with continuous kernel
$\Ks_m(x,y,\la)$ with asymptotic expansion on the diagonal
$M\subset M\times_{\pi}M$, summing over differential form degree
$p$,
\begin{eqnarray}\label{e:resolventkernelexpansion}
\Ks_m(x,x,\la) & \sim & \sum_{p=0}^{\dim M}\left(\sum_{j\geq 0,\,
[p,k]}
A_{j,[p,k]}(z,x)(-\la)^{\frac{w_k + n -j}{r}-(m+k)}\right. \nonumber \\
&+ & \left.\sum_{l\geq 0,\, [q,k]} (A^{'}_{l,[p,q]}(z,x)\log \la +
A^{''}_{l,[p,q]}(z,x))(-\la)^{-l-(m+q)}\right)
\end{eqnarray}
where $k,q\in\{0,1,\ldots,\dim M\}$ and $[p,k]
=(p_{i_1},\ldots,p_{i_k})$ an ordered multi-index of $k$
non-negative integers $p_{i_1}<\ldots < p_{i_k}$ with
$$|[p,k]| := p_{i_1}+\ldots +p_{i_k} = p
\hskip 7mm {\emph and} \hskip 7mm w_k = v_{i_1} + \ldots v_{i_k} \
,$$  and where the coefficients $A_{j,[p,k]}, A^{'}_{l,[p,k]}$ are
locally determined and $A^{''}_{l,[p,k]}$ globally determined
sections of the bundle
$$\left(\pi^*(\wedge T^{*}B)\otimes |\wedge_{\pi}|\right)\otimes\End(\Ee)
\subset \wedge^p T^{*}M\otimes\End(\Ee)$$ over $M$, depending
smoothly on  $z = \pi(x)$.

Consequently, taking the supertrace one has an asymptotic
expansion as $\la\to \o$ in $\G_{\th}$, summing over form degree
$d=p-n\geq 0$,
\begin{eqnarray}\label{e:strexpansion}
\str(\partial^{m-1}_{\la}(\Fs - \la \Is)\ii\,)(z) & \sim &
\sum_{d=0}^{\dim B}\left(\sum_{j\geq 0,\, [d,k]}
\a_{j,[d,k]}(z)(-\la)^{\frac{w_k + n -j}{r}-(m+k)}\right. \nonumber \\
& &  + \left.\sum_{l\geq 0,\, [d,q]} (\a^{'}_{l,[d,q]}(z)\log \la
+ \a^{''}_{l,[d,q]}(z))(-\la)^{-l-(m+q)}\right) \ ,
\end{eqnarray}
where the coefficients $\a_{j,[d,k]} =
\int_{M/B}\str(A_{j,[d+n,k]})$, and similarly $\a^{'}_{l,[d,q]},
\a^{''}_{l,[d,q]}$, are elements of $\Aa^{d}(B)$.

If $v_i\in \Nf$, then \eqref{e:strexpansion} can be written more
economically as
\begin{equation}\label{e:strexpansion2}
\str(\partial^{m-1}_{\la}(\Fs - \la \Is)\ii\,)(z) \sim
\sum_{d=0}^{\dim B}\left(\sum_{j\geq 0} \b_{j,d}(z)(-\la)^{\frac{w
+ n -j}{r}-m} + \sum_{l\geq 0} (\b^{'}_{l,d}(z)\log \la +
\b^{''}_{l,d}(z))(-\la)^{-l-m}\right) \
\end{equation}
$\b_{j,d}, \b^{'}_{l,[d,q]}, \b^{''}_{l,[d,q]}\in\Aa^{d}(B)$.
\end{thm}
\begin{proof}
Let $\Ws = \Fs - \Ps \in \Aa^{>0}(B,\Psi(\Ee))$. Then from
\eqref{e:resolventformula} we have
\begin{eqnarray}
(\Fs - \la \Is)\ii & = & \sum_{k=0}^{\dim B} (-1)^k(\Ps - \la
\Is)\ii
\left(\Ws(\Ps - \la \Is)\ii\right)^k  \label{e:resexpansion0} \\
& = & \sum_{k=0}^{\dim B} (-1)^k(\Ps - \la \Is)\ii
\left(\sum_{i=0}^{\dim B} \Ws_{[i]}(\Ps - \la \Is)\ii\right)^k
 \label{e:resexpansion} \ .
\end{eqnarray}
Hence
\begin{equation*}
(\Fs - \la \Is)\ii_{[d]} = (-1)^d\sum_{p_1+\ldots +p_k = d} (\Ps -
\la \Is)\ii \Ws_{[p_1]}(\Ps - \la \Is)\ii \ldots \Ws_{[p_k]}(\Ps -
\la \Is)\ii \ .
\end{equation*}
And consequently,
\begin{equation}\label{e:resolventpower}
\frac{(-1)^d}{(m_0 - 1)!\ldots (m_k - 1)!} \
\partial^{m-1}_{\la}(\Fs - \la \Is)\ii_{[d]}
\end{equation}
$$= \sum_{p_1+\ldots +p_k = d} \left(\sum_{m_0
 + \ldots +m_k = m+k}(\Ps - \la \Is)^{-m_0} \, \Ws_{[p_1]}(\Ps - \la
\Is)^{-m_1} \ldots \Ws_{[p_k]}(\Ps - \la \Is)^{-m_k}\right) \ .$$
For clarity we have assumed that with respect to any local
trivialization
\begin{equation}\label{e:onecomponentatatime}
\Ws_{[p_j]} = w_{[p_j]}\otimes Q_j
\end{equation}
with $w_{[p_j]}$ a $p_j$-form; for the general case take finite
sums of the following expansions.

The task at hand, then, is to compute the asymptotic supertrace of
the operator
\begin{equation*}
\Rs(\la,m) = (\Ps - \la \Is)^{-m_0} \, \Ws_{[p_1]}(\Ps - \la
\Is)^{-m_1} \ldots \Ws_{[p_k]}(\Ps - \la \Is)^{-m_k} \ .
\end{equation*}
First,  $\Rs(\la,m)$ is a smooth family of $\pdo$s of order
\begin{equation}\label{e:sum}
v_1 + \ldots + v_k - (m_0 +1)r - \ldots - (m_k +1)r = w_k - (m+k)r
\ .
\end{equation}
 To satisfy \eqref{e:resolventpowers} we hence need
$$ w_k - (m+k)r  < -n, \hskip 10mm k= 0, \ldots,\dim B \ ,$$
and so taking $m >(w+n)/r$ will do.

For simplicity we will treat only the case where $r$ is a positive
integer, the general case is obtained by a generalization
explained in \cite{GH}. The proof follows broadly Theorem(2.7) of
\cite{GS1}, and we use without comment the notation introduced
there  $S^{\gamma,l}$ (or $S_z^{\gamma,l} = S^{\gamma,l}(M_z)$)
for the symbol spaces (on $M_z$) \footnote{Optimally, the proof
here can be given for the vertical analogues of those symbol
spaces on a fibration of manifolds $\pi:M\to B$, but since a full
presentation of that generalized calculus is quite long we content
ourselves here with a pointwise argument.}. Let $\Bs = {\rm
OP}({\bf b})$ be the parametrix for $\Ps-\la\Is$ constructed from
\eqref{e:parametrix}. Then  $(\Ps-\la\Is)\Bs = I - \Rs$ with
pointwise $\Rs(z)\in {\rm OP}(S_z^{-\o,-r})$, while  for large
$\la\in\G_{\th}$ \eqref{e:inverse} implies
$$(\Ps-\la\Is)^{-m_i} = \Bs^{m_i} + \Bs\sum_{j\geq 1}\Rs_{(j)}$$
in $\G(B,\Psi(\Ee))$ with $\Rs_{(j)}(z)\in {\rm
OP}(S_z^{-\o,-(j+1)r})$. Consequently
\begin{equation}\label{e:weaklypolyhomogeneous}
(\Ps - \la \Is)^{-m_0} \, \Ws_{[p_1]}(\Ps - \la \Is)^{-m_1} \ldots
\Ws_{[p_k]}(\Ps - \la \Is)^{-m_k}
\end{equation}
$$=
\Bs^{m_0}\Ws_{[p_1]}\Bs^{m_1}\Ws_{[p_2]}\Bs^{m_2}\ldots
\Ws_{[p_k]}\Bs^{m_k} + \sum_{j\geq
1}T(\Bs,\Ws_{[p_1]},\ldots,\Ws_{[p_k]},\Rs_{(j)})$$ is a smooth
family of weakly polyhomogeneous $\pdo$s with differential form
degree $d$. The vertical operator coefficient of
$T(\Bs,\Ws_{[p_1]},\ldots,\Ws_{[p_k]},\Rs_{(j)})(z)$  is in ${\rm
OP}(S_z^{-\o,-(j+1)r})$, while by Theorem (1.12) of \cite{GS1} the
expansion of the terms of local vertical symbol of
$T(\Bs,\Ws_{[p_1]},\ldots,\Ws_{[p_k]},\Rs_{(j)})$ are expansions
in integer powers $\la = \mu^r$. By Theorem (2.1) of \cite{GS1}
the kernel
\begin{equation*}
K_{T(\Bs,\Ws_{[p_1]},\ldots,\Ws_{[p_k]},\Rs_{(j)})} \in
\G(M\times_{\pi}M, \pi^*(\wedge T^*B)\otimes(\Ee\otimes
|\Lambda_{\pi}|^{1/2})\otimes(\Ee^*\otimes |\Lambda_{\pi}|^{1/2}))
\end{equation*}
has an expansion as $\la\to\o$ in $\G_{\th}$
\begin{equation}\label{e:KT}
K_{T(\Bs,\Ws_{[p_1]},\ldots,\Ws_{[p_k]},\Rs_{(j)})}(z,x,x,\la)
\sim \sum_{\s\geq 0}C_{j,\s}(z,x)(-\la)^{-j-(m+k)-\s} \ ,
\end{equation}
with $C_{j,\s}\in\G(M,\left(\pi^*(\wedge T^{*}B)\otimes
|\wedge_{\pi}|\right)\otimes\End(\Ee))$ and highest power
$(-\la)^{-1-m-k}$.

On the other hand, with the assumption
\eqref{e:onecomponentatatime}, in a local trivialization let
$\s(\As)$ denote the local vertical symbol of the  $\pdo$ family
in $\G(B,\Psi(\Ee))$ coefficient to $\As\in\Aa^d(B,\Psi(\Ee))$, we
have $$\s((\Ps - \mu^r)^{-m_j})_z\in S_z^{-rm_j,0}\cap
S_z^{0,-rm_j} \ .$$ Hence with $q  = \s((\Ps - \la \Is)^{-m_0} \,
\Ws_{[p_1]}(\Ps - \la \Is)^{-m_1} \ldots \Ws_{[p_k]}(\Ps - \la
\Is)^{-m_k})$, \eqref {e:sum} and the symbol calculus of
\cite{GS1} imply that
$$q_z\in S_z^{w_k -r(m+k),0}\cap S_z^{w_k,-r(m+k)} \ $$
with an expansion
$$q(z,x,\xi,\mu) \sim \sum_{j\geq 0}q_{w_k -r(m+k)-j}(z,x,\xi,\mu) \ $$
where $q_{w_k -r(m+k)-j}\in S_z^{w_k -j,(m+k)r}$. From \cite{GS1}
Thm (1.12), a symbol $p\in S^{k,d}$ has a Taylor expansion as
$\mu\to\o$
$$p(x,\xi,\mu) = \sum_{j=0}^{J}p^{(j)}(x,\xi) \mu^{d-j} +
O((1+|\xi|^2)^{(k+J)/2)}\,\mu^{d-J}) \ $$ with $p^{(j)}\in
S^{k+j}$. Consequently, locally since the kernel of
$\Bs^{m_0}\Ws_{[p_1]}\Bs^{m_1}\Ws_{[p_2]}\Bs^{m_2}\ldots
\Ws_{[p_k]}\Bs^{m_k}$ restricted to the diagonal is
$$K_{{\rm OP}(q)}(z,x,x,\mu) =
\frac{i}{(2\pi)^n}\int_{\Rf^n}q(z,x,\xi,\mu) d\xi \otimes
\omega_{[d]}(z)\otimes \upsilon_z$$ where $\upsilon_z$ is a local
volume form on $U_z\subset M_z$ and $\omega_{[d]}$ the local
coefficient $d$-form on the base $B$, and $\mu = \la^{1/r}$
relative to $R_{\th}$, then by splitting the integral into three
summands for $|\xi|\geq |\mu|, |\xi|\leq 1$ and $1\leq|\xi|\leq
|\mu|$ we obtain by the proof of \cite{GS1} Theorem (2.1) a kernel
expansion
\begin{equation}\label{e:qkernelexpansion}
K_{\Bs^{m_0}\Ws_{[p_1]}\Bs^{m_1}\ldots
\Ws_{[p_k]}\Bs^{m_k}}(z,x,x,\mu)
\end{equation}
$$\sim \sum_{j\geq 0} B_j (z,x)(-\la)^{\frac{w_k + n -j}{r}-(m+k)}
+ \sum_{l\geq 0} (B^{'}_l (z,x)\log \la + B^{''}_l
(z,x))(-\la)^{-l-(m+q)} \ ,$$ where $B_{j}, B^{'}_l, B^{''}_l
\in\G(M,\left(\pi^*(\wedge T^{*}B)\otimes
|\wedge_{\pi}|\right)\otimes\End(\Ee))$.

From \eqref{e:resolventpower}, \eqref{e:weaklypolyhomogeneous},
\eqref{e:KT} and \eqref{e:qkernelexpansion} we obtain the
expansion \eqref{e:resolventkernelexpansion}, from which the
remaining statements are immediate consequences.
\end{proof}

\noi {\bf Remark.}  \newline \noi {\rm [1]} With obvious
modifications to the powers of $\la$, for an auxiliary $\As\in
\Aa(B,\Psi(\Ee))$  the resolvent supertrace expansion
\eqref{e:strexpansion} extends to
\begin{eqnarray}\label{e:Astrexpansion}
\str(\partial^{m-1}_{\la}(\As(\Fs - \la \Is)\ii\,) \ .
\end{eqnarray}
Further, generalizations yield expansions for families of
pseudodifferential operators with powers of $\log|\xi|$ in the
homogeneous terms of the vertical symbol, which is essential, for
example, for zeta determinant form formulae and higher
multiplicative anomalies.\vskip 2mm

\vskip 3mm

It is important to see how the asymptotic expansions transform
with respect to the rescaling by $t>0$ of \cite{B,BGV}
\begin{equation*}
\d_t :\Aa(B,\pi_*(\Ee)) \to  \Aa(B,\pi_*(\Ee)), \hskip 10mm \d_t
\omega_{[i]} = t^{-i/2}\omega_{[i]} \  .
\end{equation*}
$\d_t$ induces an automorphism of $\Aa(B,\Psi(\Ee))$ given by $
\d_t(\As) = \d_t\cdot\As\cdot\d_t\ii.$ Let
$\Fs\in\Aa(B,\Psi(\Ee))$ satisfying the assumptions of
\thmref{t:0}, and  define
\begin{equation}\label{e:rescaledF}
\Fs_t = t\d_t(\Fs) \ .
\end{equation}
Then
$$
\str(\partial^{m}_{\la}(\Fs_t - \la \Is)\ii\,) =
\d_t(\str(\partial^{m}_{\la}(t\Fs - \la \Is)\ii\,)) =
t^{-m-1}\d_t(\str(\partial^{m}_{\la}(\Fs - \la t\ii \Is)\ii\,)) \
,$$ and  from \eqref{e:strexpansion}, since the coefficients are
in $\Aa^d(B)$, we obtain
\begin{eqnarray}\label{e:strrescaled}
\str(\partial^{m-1}_{\la}(\Fs_t - \la \Is)\ii\,) & \sim &
\sum_{d=0}^{\dim B}\left(\sum_{j\geq 0,\, [d,k]}
\a_{j,[d,k]}(-\la)^{\frac{w_k + n -j}{r}-(m+k)} \
t^{\frac{j - w_k - n}{r}+k -1 -\frac{d}{2}})\right. \nonumber \\
& &  + \left.\sum_{l\geq 0,\, [d,q]} (\a^{'}_{l,[d,q]}\log (\la
t\ii) + \a^{''}_{l,[d,q]})(-\la)^{-l-(m+q)} \
t^{l+q-1-\frac{d}{2}}\right) \ ,
\end{eqnarray}
while in the case $v_i\in \Nf$, then \eqref{e:strexpansion2}
rescales to
\begin{eqnarray}\label{e:strrescaled2}
\str(\partial^{m-1}_{\la}(\Fs_t - \la \Is)\ii\,) & \sim &
\sum_{d=0}^{\dim B}\left(\sum_{j\geq 0} \b_{j,d}(-\la)^{\frac{w +
n -j}{r}-m} \ t^{\frac{j - w - n}{r} -1 -\frac{d}{2}} \right.
\nonumber \\ & & + \left. \sum_{l\geq 0} (\b^{'}_{l,d}\log \la
(t\ii) + \b^{''}_{l,d})(-\la)^{-l-m} \ t^{l-1-\frac{d}{2}}\right)
\ .
\end{eqnarray}

\section{Zeta Forms and Zeta Determinant Forms}

Let $\Fs\in \Aa(B,\Psi(\Ee))$ be a smooth family of elliptic
$\pdo$s of constant order $(r>0,v_1,v_2,\ldots,v_{\dim B})$ with
parameter $\la\in\G_{\th}$. Then, with $\Ps = \Fs_{[0]}$, for
$\la\in R_{\th}$ sufficiently large $\Ps-\la \Is$ is a form degree
zero family of invertible $\pdo$s of positive order $r$. If
$\Ps-\la\Is$, and hence $\Fs-\la\Is$, is invertible for all $\la$
in the spectral cut $R_{\th} = R_{\th}\backslash\{0\}$ then the
angle $\th$ is called an {\em Agmon angle} for $\Fs$
%\footnote{In
%the case of a single operator with spectral cut $R_{\th}$, an
%Agmon angle can be found by a small perturbation of $\th$, but for
%a family of elliptic $\pdo$s the existence must be postulated. The
%existence of an Agmon angle is not actually a necessary
%requirement for the constructions here, but we include it to avoid
%some somewhat arid technical analysis.}

We assume for the moment that
\begin{equation}\label{e:orderlessthanr}
 v_k = \ord(\Fs_{[k]}) \leq r, \hskip 10mm k=1,\ldots,\dim B \ .
\end{equation}
With \eqref{e:orderlessthanr} and using the expansion
\eqref{e:resexpansion} we obtain an operator norm estimate in
$\Aa(B)$ as $\la\to\infty$ in $\G_{\th}$
\begin{equation}\label{e:Olambdainverse}
\| (\Fs - \la \Is)\ii\|^{(l)}_{M/Z} = O(|\la|\ii) \ ,
\end{equation}
where for $l\in\Rf$
\begin{equation}\label{e:verticaloperatornorm}
\| \ . \ \|^{(l)}_{M/B} : \Aa(M,\Psi^{0}(\Ee)) \too \Aa(B)
\end{equation}
is the vertical operator Sobolev norm associated to the vertical
metric
\begin{equation}\label{e:verticalmetric}
| \ . \ |_{M/B} : \Aa(B,\pi_*(\Ee)) \too \Aa(B) \ , \hskip 10mm |
\pi^*(\a)\otimes\Psi\otimes\upsilon |_{M/B} = \a
\int_{M/B}|\psi|^2\,\upsilon^2 \ ,
\end{equation}
defined independently of the representation of a section as tensor
product
$$\psi\otimes\pi^*(\a)\otimes\upsilon\in \G(M,\pi^*(\wedge T^*
B)\otimes\Ee\otimes |\wedge|_{\pi}) \ . $$ Pointwise for $z\in B$
the metric \eqref{e:verticalmetric} reduces on the fibre
$\G(M_z,\Ee\otimes|\wedge_{M_z}|)$ to the canonical metric
$|\psi_z|^2 = \int_{M_z}|\psi_z(x)|^2.$

On the right side of \eqref{e:Olambdainverse}, for $a:\Cf\to
\Aa(B)$ we write $a(\la) = O(f(\la))$ if for each $l\in\Nf$ and
relatively compact subset $U$ of $B$, there is a constant $C(l,U)$
such that $\|a(\la)\|_l \leq C(l,U)f(\la)$, where $\|\ . \ \|_l$
is the $C^l$ norm. The proof of \eqref{e:Olambdainverse} follows
that of the classical result for a single operator \cite{S,Sh}
using the vertical parametrix $\Bs$ in \eqref{e:inverse}. As for
Seeley's analysis of the single operator case \cite{S}, it follows
if $\Fs$ has Agmon angle $\th$, that for $\re(s)>0$ the complex
powers can be defined by
\begin{equation}\label{e:powers}
\Fs_{\th}\si = \frac{i}{2\pi}\int_{C}\la_{\th}\si\, (\Fs -
\la\Is)\ii \, d\la \ ,
\end{equation}
where $\la_{\th}$ is the branch of $\la\si$ defined by $
\la_{\th}\si = |\la|\si e^{-is\arg(\la)}$ with $\th - 2\pi \leq
\arg(\la) < \th$, and where $C$ is the negatively oriented contour
 which is the boundary of a sector
\begin{equation*}%\label{e:sector}
\Lambda_{\th,\d} = \{z\in \C \ | \ |\arg(z)-\th| \leq \d \;\;{\rm
or}\;\; |z| \leq \rho\}
\end{equation*}
with $\d$ so chosen that $\Lambda_{\th,\d}$ contains no
eigenvalues of the operators $\Ps_z$ and $\rho$ such that $(\Fs -
\la\Is)\ii$ is defined and holomorphic for $ 0 < |\la| < \rho +
\e$ for some $\e>0$.

The estimate \eqref{e:Olambdainverse} shows that $\Fs\si$
converges in each vertical Sobolev norm and hence defines an
operator from $\Aa(B,\pi_*(\Ee))$ into $\Aa(B,\pi_*(\Ee))$. From
\eqref{e:resexpansion}  we have for $\re(s)>0$
\begin{equation}\label{e:powersexpansion}
\Fs_{\th}\si = \sum_{\stackrel{k}{p_1+\ldots +p_m = k}}
\frac{i}{2\pi}\int_{C}\la_{\th}\si\, (-1)^k (\Ps - \la \Is)\ii
\Ws_{[p_1]}(\Ps - \la \Is)\ii \ldots \Ws_{[p_m]}(\Ps - \la
\Is)\ii\, d\la .
\end{equation}
Each of the summands in \eqref{e:powersexpansion} is smooth family
of $\pdo$s, of differential form degree $k$, represented locally
by a sum of vertical polyhomogeneous symbols
\begin{equation}\label{e:powersymbols}
\frac{i}{2\pi}\int_{C}\la_{\th}\si\, {\bf
b}_{m_1}[\la](z,x,\xi)\circ {\bf w}_{p_l,i_1}(z,x,\xi)\ldots {\bf
w}_{p_m,i_m}(z,x,\xi)\circ{\bf b}_{m_j}[\la](z,x,\xi)\, d\la
\end{equation}
with $\Bs = {\rm OP}(\sum {\bf b}_i[\la])$ the parametrix for
$\Ps-\la\Is$ and $\Ws_{p_i} = {\rm OP}(\sum_{\s}{\bf w}_{p_i,\s})$
where ${\bf w}_{p_i,\s}(z,x,\xi)$ is homogeneous in $\xi$ of order
$\nu_i-\s$ and ${\bf b}_i[\la])(z,x,\xi)$ homogeneous in
$(\xi,\la^{1/r})$ of degree $r-i$, i.e. ${\bf
b}_i[t^r\la])(z,x,t\xi) = t^{r-i}{\bf b}_i[\la])(z,x,\xi)$ for
$t>0$, $\la,t^r\la\in \Lambda_{\th}$. The degree of homogeneity of
\eqref{e:powersymbols} is computed by replacing $\xi$ by $t\xi$
and $\la$ by $t^r\mu$ in the integrand of symbol products. In
particular, setting $m_1 = \ldots = m_j = 0$, $\s = 0$, we have
that the principal symbol has degree $\nu_{i_1} + \ldots +
\nu_{i_m} - (s+k)r$.

Proceeding in this way, applying the standard methods of
\cite{S,Sh} and the remark following \propref{p:resolventformula},
we obtain the following fact.
\begin{lem}\label{lem:powersarepdos}
The vertical complex power defined for $\re(s)>0$  by
\eqref{e:powers} (and generally without assumption
\eqref{e:orderlessthanr} for $\re(s) >>0$, see below,) is a smooth
family of $\pdo$s of mixed differential form degree
$\Fs\si\in\Aa(B,\Psi(\Ee))$ such that if ${\bf f}[\la] \sim
\sum_{j}{\bf f}[\la]_j$ is a local vertical polyhomogeneous symbol
representing $(\Fs-\la\Is)\ii$, then ${\bf f}_{\th}\si \sim
\sum_{j}{\bf f}\si_{\th,j}$ represents $\Fs\si$, where
$${\bf f}\si_{\th,j}(z,x,\xi) =
\frac{i}{2\pi}\int_{C}\la_{\th}\si\, {\bf f}[\la]_j(z,x,\xi) \,
d\la \ .
$$
\end{lem}

\vskip 2mm

The zeta form for $\Fs$ can now be constructed as follows. Since
\eqref{e:Olambdainverse} implies for $\re(s)>0$ the operator norm
estimate in $\Aa(B)$
\begin{equation*}
\| \la^{m-s}\, \dd_{\la}^{m-1}(\Fs - \la \Is)\ii\|^{(l)}_{M/Z} =
O(|\la|\ii) \ ,
\end{equation*}
as $\la\to\infty$ along $C$, we can integrate by parts in
\eqref{e:powers} to obtain
\begin{equation}\label{e:F-by-parts}
  \Fs\si  =  \frac{1}{(s-1)\ldots (s-m)}.\frac{i}{2\pi}\int_{C}\la^{m-s}\,
  \dd_{\la}^{m}(\Fs - \la\Is)\ii \, d\la \ .
\end{equation}
Since
\begin{equation*}
\dd_{\la}^{m}(\Fs - \la\Is)\ii = \sum_{\stackrel{k=0}{m_0 + \ldots
+ m_k = m}}^{\dim B} (-1)^k\dd_{\la}^{m_0}((\Ps - \la \Is)\ii)
\Ws\,\dd_{\la}^{m_2}((\Ps - \la \Is)\ii) \ldots \Ws \,
\dd_{\la}^{m_k}((\Ps - \la \Is)\ii)\,d\la
\end{equation*}
and $\dd_{\la}^{m_i}(\Ps - \la \Is)\ii \in
\Aa(B,\Psi^{-rm_i-r}(\Ee))$, then by taking $m\geq N$ for
sufficiently large $N$ we may ensure an estimate $\|
\dd_{\la}^{m}(\Fs - \la \Is)\ii\|^{(l)}_{M/Z} = O(|\la|\ii) $
without assuming \eqref{e:orderlessthanr}. For the general case,
we hence define $\Fs\si$ for $\re(s) > m \geq N$ by
\eqref{e:F-by-parts}.

Moreover, \eqref{e:F-by-parts} and \thmref{t:0}, equation
\eqref{e:resolventpowers}, show for $\re(s) > \frac{w+n}{r}$ that
$\Fs\si \in \Aa(B,\Psi^{< -n}(\Ee))$ is a smooth family of trace
class $\pdo$s with kernel $K(\Fs\si)$ continuous over the diagonal
$M\subset M\times_{\pi}M$. In that half-plane $\Fs$ therefore has
a super-zeta form
\begin{equation}\label{e:superzetaform}
\z_{\th}(\Fs,s) :=  \str(\Fs_{\th}\si) = \int_{M/B}
\str(K(\Fs_{\th}\si)(x,x)) \in \Aa(B) \ .
\end{equation}
From \eqref{e:F-by-parts} for $\re(s) > m > (w+n)/r$ we have
\begin{equation}\label{e:superzetaform2}
\z_{\th}(\Fs,s) =\frac{1}{(s-1)\ldots
(s-m)}.\frac{i}{2\pi}\int_{C}\la^{m-s}\,
\str\left(\dd_{\la}^{m}(\Fs - \la\Is)\ii\right) \, d\la \ .
\end{equation}
We can use \eqref{e:superzetaform2} to write down the singularity
structure of the meromorphic continuation of the zeta form to all
$s\in \Cf$. To do so requires the assumption that $\Ps =
\Fs_{[0]}$ is smooth family of $\pdo$s such that $\Ker(P_z)$ has
constant kernel dimension. Consequently the meromorphically
continued zeta form $\z_{\th}(\Fs,s)|^{{\rm mer}}$ is only defined
for families of $\pdo$s with Agmon angle in $\Aa(B,\Psi(\Ee))$
modulo the regularizing subalgebra $\Aa(B,\Psi^{-\o}(\Ee))$. This
is essentially equivalent to $\z_{\th}(\Fs,s)|^{{\rm mer}}$ being
a characteristic class map on K-theory, taking values in $H^*(B)$.

Thus we assume that the family of $\pdo$ projections onto the
kernels define a smooth family of smoothing operators
$\Pi\in\Aa(B,\Psi^{-\o}(\Ee))$, defining a smooth finite-rank
superbundle $\Ker(\Ps)$ over $B$. With this assumption,
\eqref{e:resexpansion0} implies that at $\la =0$ the resolvent
trace form is meromorphic with Laurent expansion
\begin{equation}\label{e:Laurentatzero}
\str\left(\dd_{\la}^m(\Fs - \la \Is)\ii\right)  = \sum_{k=0}^{\dim
B} (-1)^k \frac{(m+k)!}{k!}(- \la)^{-k-1-m}
\str\left((\Pi\cdot\Ws\cdot\Pi)^k\right) + O(|\la|^m) \ .
\end{equation}

The asymptotic expansion \eqref{e:strexpansion} as $\la\to\o$ in
$\G_{\th}$ along with the expansion \eqref{e:Laurentatzero} at
$\la =0$ now imply by a standard transition argument, for example
\cite{GS2} Proposition (2.9), that $\z_{\th}(\Fs,s)$ extends
meromorphically to $\Cf$ with the singularity structure

\vskip 1mm

\begin{equation}\label{e:singularitystructure}
 \frac{\pi}{\sin(\pi s)}\z_{\th}(\Fs,s)|^{{\rm mer}} \sim
 \end{equation}
 $$ -\sum_{j=-\dim B -1}^{-1}
  \frac{\str\left((\Pi\cdot\Ws\cdot\Pi)^{-j-1}\right)}{(s - j -
  1)} \ \ + $$
 $$ \sum_{d=0}^{\dim B}\left(\sum_{j\geq 0,\,[d,k]}
\frac{a_{j,[d,k]}}{\left(s + k + \frac{j - n - w_k}{r} -1
\right)}\right.   +  \left. \sum_{l\geq 0,\,
[d,q]}\frac{a^{'}_{l,[d,q]}}{(s + l + q - 1)^2}  +
\frac{a^{''}_{l,[d,q]}}{(s + l + q - 1)}\right) \ ,$$ \vskip 4mm
\noi with coefficients $a_{j,[d,k]}, a^{'}_{l,[d,q]},
a^{''}_{l,[d,q]}\in \Aa^{d}(B)$ related to the coefficients of
\eqref{e:strexpansion} by universal multiplicative constants.
Specifically,
\begin{equation}\label{e:singcoeffs}
 a_{j,[d,k]} =  \G(\frac{j - n - w_k}{r} + k)\G(\frac{j - n - w_k}{r} + k + m )\ii
 \, \a_{j,[d,k]}(m)
\end{equation}
independently of m, with $\G(s)$ the Gamma function. In the more
general case where we allow non-constant order $\nu_k\in\G(B,\Rf)$
then the factors are replaced by the corresponding universal
functions.

If $v_i\in \Nf$, then \eqref{e:singularitystructure2} takes the
simpler form

\vskip 1mm

\begin{equation}\label{e:singularitystructure2}
 \frac{\pi}{\sin(\pi s)}\z_{\th}(\Fs,s)|^{{\rm mer}} \sim
 \end{equation}
 $$-\sum_{j=-\dim B -1}^{-1}
  \frac{\str\left((\Pi\cdot\Ws\cdot\Pi)^{-j-1}\right)}{(s - j -
  1)^{k+1}} \ \ + $$
 $$\sum_{d=0}^{\dim B}\left(\sum_{j\geq 0}
\frac{b_{j,d}}{\left(s + \frac{j - n - w}{r} -1 \right)}\right. +
\left. \sum_{l\geq 0}\frac{b^{'}_{l,d}}{(s + l - 1)^2} +
\frac{b^{''}_{l,d}}{(s + l - 1)}\right) \ ,$$ \vskip 4mm \noi with
coefficients $b_{j,d}, b^{'}_{l,d}, b^{''}_{l,d}\in \Aa^{d}(B)$
related to the $\b_{j,d}, \b^{'}_{l,d}, \b^{''}_{l,d}$ by
constants, and
\begin{equation}\label{e:singcoeffs2}
 b_{j,d} =  \G(\frac{j - n - w}{r})\G(\frac{j - n - w}{r} + m )\ii
 \, \b_{j,d}(m)
\end{equation}

The pole structure \eqref{e:singularitystructure} can also be
computed directly from the meromorphically continued symbol
representation of $\Fs\si$ in \lemref{lem:powersarepdos}.

\begin{defn}
A family of $\pdo$s $\Fs\in\Aa(B,\Psi(\Ee))$ admitting an Agmon
angle $\th$ is said to be $\z$-admissible if when $l+q-1\in
\{0,1,\ldots,\dim B\}$ then $a^{'}_{l,[d,q]} = 0$ for $d\geq 1$ in
\eqref{e:singularitystructure} (for $d=0$ this is guaranteed by
the ellipticity of $P = \Fs_{[0]}$). Similarly, for
\eqref{e:singularitystructure2} this requires $b^{'}_{l,d} = 0$
for $l-1 \in \{0,1,\ldots,\dim B\}$.
\end{defn}
This ensures that $\z_{\th}(\Fs,s)|^{{\rm mer}}$ is holomorphic
for $s$ around $0, 1,\ldots, \dim B$. This property is needed for
the differential $\zeta$ form, but when $\Fs$ is the curvature of
a superconnection is irrelevant at the cohomological level, since
in that case the forms are all exact.

The complex powers $\Fs\si$ defined by \eqref{e:powers} for
$\re(s) > 0$ if \eqref{e:orderlessthanr} holds, and  in general by
\eqref{e:F-by-parts} for $\re(s) > m$ if it does not, are extended
by Seeley's method \cite{S} to all $s\in\Cf$ by choosing any
positive integer $N$ with $\re(s) + N > m$ and defining
\begin{equation}\label{e:genpowers}
    \Fs\si = \Fs^{-s-N}\Fs^N \in \Aa(B,\Psi(\Ee)) \ .
\end{equation}
More precisely, the map $s\mtoo K(\Fs\si)$ assigning to $\Fs\si$
its (distributional) kernel is a holomorphic map of $\{s \ | \
\re(s) > (w+n)/r\}$ into (in a local trivialization) matrices of
continuous functions. Restricted to any compact subset $V$ of
$M\times M\backslash \Delta(M)$ the map $s\in\Cf\mtoo
K(\Fs\si)_{|V}$ is holomorphic map from $\Cf$ to smooth matrices,
while along the diagonal $s\mtoo K(\Fs\si)(x,x)$ is a meromorphic
function on all of $\Cf$ with discrete poles at the points
indicated in \eqref{e:singularitystructure}.

We then define the logarithm of $\Fs$ to be the smooth vertical
family of {\em log-polyhomogeneous} $\pdo$s
\begin{equation*}
\log_{\th} \Fs := - \dd_s |_{s=0}\Fs_{\th}\si \in \Aa(B,\Psi_{{\rm
log}}(\Ee))\ .
\end{equation*}
Thus, omitting the $\th$ subscript, $\dd_s \Fs\si = -
\log\Fs\,\Fs\si$, where for $\re(s) > 0$ if
\eqref{e:orderlessthanr} holds
\begin{equation*}
 \log\Fs\,\Fs\si = \frac{i}{2\pi}\int_{C}\log\la \,\la\si\, (\Fs -
\la\Is)\ii \, d\la \ ,
\end{equation*}
and similarly using \eqref{e:F-by-parts} for the general case.

Here $\Aa(B,\Psi_{{\rm log}}(\Ee))$ is the extension of
$\Aa(B,\Psi(\Ee))$ to operators represented by vertical
log-polyhomogeneous symbols. This means that with respect to local
coordinates on  $\Ee$, an operator $\Ts\in\Aa(B,\Psi_{{\rm
log}}(\Ee))$ is represented by a vertical symbol ${\bf
t}\in\G\left( \, (U_M\times_{\pi} U_M) \times
\Rf^{n}\backslash\{0\}, \ \pi^*(\wedge T^*
U_B)\otimes\Rf^N\times(\Rf^N)^*\right)$ of the form
\begin{equation}\label{e:verticallogsymbol}
{\bf t}(z,x,\xi) \sim\sum_{j\geq 0}\sum_{p=0}^1 {\bf
t}_{-j,p}(z,x,\xi)\log^p|\xi| \ ,
\end{equation}
with ${\bf t}_{-j,p}(z,x,\xi)$ a vertical homogeneous symbol.

It is readily verified that $\log_{\th}\Fs$ is log-polyhomogeneous
locally represented by the vertical log-polyhomogeneous symbol
$\log {\bf f} \sim\sum_{j\geq 0}\log_j{\bf f}$ with
$$\log_j{\bf f}(x,\xi) = \frac{i}{2\pi}\int_{C}\log \la\,
{\bf b}[\la]_j(z,x,\xi) \, d\la \ $$ and $ \dd_s {\bf f}\si =
(\log {\bf f})\circ {\bf f}\si$ and furthermore $(\log\Fs)_{[0]} =
\log\Ps$.

\begin{defn}
For $\z$-admissible $\Fs\in\Aa(B,\Psi(\Ee))$ with Agmon angle
$\th$, the zeta-determinant form $\det_{\z,\th}\Fs\in\Aa(B)$ is
defined by
\begin{equation}\label{e:trlog}
\log{\rm det}_{\z,\th}\Fs\in\Aa(B)  = - \dd_s|^{{\rm mer}}_{s=0}\
\str(\Fs\si)  =  \str(\log\Fs \, . \, \Fs\si)|^{{\rm mer}}_{s=0}
  \ .
\end{equation}
\end{defn}

\vskip 2mm

\lemref{p:zetaperturbation} is an immediate corollary of the
following.
\begin{lem} One has
 $ \z(\Fs,s)_{[0]} = \z(\Ps,s).$
\end{lem}

\noi To see this, notice from \eqref{e:F-by-parts}  and
\eqref{e:resexpansion0} that for $\re(s) >>0$
$$
\Fs\si = \Ps\si + \sum_{k=0}^{\dim B}\frac{1}{(s-1)\ldots
(s-m)}.\frac{i}{2\pi}\int_{C}\la^{m-s}\, \dd_{\la}^{m}\left(\,(\Ps
- \la\Is)\ii (\Ws (\Ps - \la\Is)\ii )^k\,\right) \, d\la \ , $$
and hence that for $m > (w + n)/r$
$$\zeta(\Fs,s)|^{{\rm mer}}  = \z(\Ps,s)|^{{\rm mer}} \ + $$ $$
\sum_{k=0}^{\dim B}\frac{1}{(s-1)\ldots
(s-m)}.\frac{i}{2\pi}\int_{C}\la^{m-s}\,
\str\left(\dd_{\la}^{m}(\,(\Ps - \la\Is)\ii (\Ws (\Ps - \la\Is)\ii
)^k\,)\right) \, d\la \left. \right|^{mer}\ .$$ The second term on
the right is meromorphic on $\Cf$, and holomorphic at zero, and of
non-zero form degree, since $\Ws\in\Aa^{>0}(B,\Psi(\Ee))$, with
the pole structure of \eqref{e:singularitystructure} but with
$d\geq 1$.

\vskip 1mm

Replacing $\Fs$ by $\Fs_t = t\d_t(\Fs)$ in
\eqref{e:singularitystructure}, we can use \eqref{e:strrescaled}
to write down the $t$-rescaled singularity structure. This means
that the rescaled left side of \eqref{e:singularitystructure}
minus the rescaled sums on the right-side for $j\leq N$ terms is
$O(t^{\frac{N - w_k - n}{r}+k -\frac{d}{2}})$. Taking the $s$
derivative and evaluating at $s=0$, this implies the following.
\begin{prop}
There is an asymptotic expansion as $t\too 0+$ in $\Aa(B)$
\begin{equation}\label{e:logdetasymptotics}
\log{\rm det}_{\z}\Fs_t \sim \sum_{d=0}^{\dim B}\sum_{j\geq 0,\,
[d,k]} c_{j,[d,k]} \, t^{\frac{j - w_k - n}{r}+k -1 -\frac{d}{2}}
+ \sum_{l\geq 0,\, [d,q]} ( c^{'}_{l,[d,q]}\log t +
c^{''}_{l,[d,q]})\,  t^{l-1-\frac{d}{2}} \ ,
\end{equation}
where the degree $d$ differential forms $c_{j,[d,k]},
c^{'}_{l,[d,q]}$ are determined locally, while $c^{''}_{l,[d,q]}$
are globally determined.
\end{prop}

\vskip 2mm

In the case where eigenvalues of the principal symbol $ {\bf
p}_{0}\in \G(T(M/B)\backslash\{0\},p^*(\End(\Ee)))$ of $\Ps$ lie
pointwise for $(x,\xi)\in T(M/B)\backslash\{0\}$ in a subsector of
the right-half plane, then, the zeta form and zeta determinant
form can be equivalently formulated by a vertical heat trace form
\begin{equation}\label{e:heattrace2}
\str(e^{-\Fs}) = \frac{i}{2\pi}\int_{\Cc}
e^{-\la}\str(\dd_{\la}^m(\Fs - \la\Is)\ii) \
  d\la \ .
\end{equation}
where $m > (w+n)/r$ and $\Cc$ is a contour coming in on a ray with
argument in $(0,\pi/2)$, encircling the origin, and leaving on a
ray with argument in $(-\pi/2,0)$. We hence obtain an asymptotic
expansion, which for brevity we state only for the case $v_i\in
\Nf$, that as $t\too 0+$
\begin{equation}\label{e:heatraceexpansion2}
\str(e^{-\Fs_t})  \sim  \sum_{d=0}^{\dim B}\left(\sum_{j\geq 0}
 \tilde{b}_{j,d} \, t^{\frac{j - w - n}{r} -1 -\frac{d}{2}} + \sum_{l\geq
0} (\tilde{b}^{'}_{l,d}\log t + \tilde{b}^{''}_{l,d}) \,
t^{l-1-\frac{d}{2}}\right) \
\end{equation}
with $\tilde{b}_{j,d}, \tilde{b}^{'}_{l,q},
\tilde{b}^{''}_{l,d}\in \Aa^{d}(B)$  related to those in the
expansion \eqref{e:singularitystructure2} by
\begin{equation}\label{e:heatzetafactors}
\tilde{b}_{j,d} = \G(\frac{j - n - w_k}{r})\ii b_{j,d} \ , \ \ \
\tilde{b}^{'}_{l,d} = \G(l)\ii b^{'}_{l,d} \ , \ \ \ \
 \tilde{b}^{''}_{l,d} = \G(l)\ii b^{''}_{l,d} \ .
\end{equation}

\section{Homotopy Properties}

A superconnection \cite{Q,B,BGV} on $\pi_*(\Ee)$ adapted to a
smooth family of formally self-adjoint elliptic $\pdo$s $\Ps =
\begin{bmatrix}
  0 & \Ps^- \\
  \Ps^+ & 0 \\
\end{bmatrix} \in\Aa^0(B,\Psi^r(\Ee))$ of order $r>0$ is a classical
$\pdo$ $\Af$ on $\Aa(B,\pi_*(\Ee))  = \G(M,\pi^*(\wedge T^*
B)\otimes \Ee\otimes|\wedge_{\pi}|^{1/2})$ of odd-parity  with
respect to the $\Zf_2$-grading, such that
\begin{equation}\label{e:superconnection}
    \Af( \omega\,\psi) = d\omega \, \psi +
    (-1)^{|\omega|}\omega\,\Af(\psi) \ ,
\end{equation}
for $\omega\in\Aa(B)$ and $\psi\in \Aa(B,\pi_*(\Ee))$, and such
that $\Af_{[0]} = \Ps \ ,$ where $\Af = \sum_{i=0}^{\dim
B}\Af_{[i]}$ and $\Af_{[i]}: \Aa^d(B,\pi_*(\Ee))\to
\Aa^{d+i}(B,\pi_*(\Ee))$ is  the component which raises form
degree by $i$. It follows from \eqref{e:superconnection} that
$\Af_{[1]}$ is a connection in the classical (ungraded) sense,
while each remaining term is a smooth family of $\pdo$s
$\Af_{[i]}\in\Aa^i(B,\Psi(\Ee))$ if $i\neq 1.$ In a local weak
trivialization $\Aa(U,\pi_*(\Ee)_{|U}) \cong \Aa(U)\otimes
\G(M_{z_0},\Ee^{z_0})$ for $z_0\in U$, $\Af$ takes the local
coordinate form $\Af_{|U} = d_{U} + \sum_I P_I dz_I,$ where $dz_I
= dz_{i_1}\ldots dz_{i_m}$ and $P^z_I$ a classical $\pdo$ on
$\G(M_{z_0},\Ee^{z_0})$.

The curvature of $\Af$ is the smooth family of $\pdo$s with
differential form coefficients $\Af^2 \in \Aa(B,\pi_*(\Ee)).$
Since $\Af^2_{[0]} = \Ps^2 =  \Ps^-\Ps^+\oplus \Ps^+\Ps^-\in
\Aa^0(B,\pi_*(\Ee))$ we have that $\Af^2$ is elliptic with Agmon
angle $\pi$.
\begin{prop}\label{p:resformclosed}
Let $m\in\Nf$ with $m > \frac{w+n}{r}$ with $w =\sum_i
\ord((\Af^2)_{[i]})$. Then for $\la\in\G_{\pi}$ the resolvent
trace form $\str(\dd_{\la}^{m-1}(\Af^2 - \la\Is)\ii)\in\Aa(B)$ is
a closed differential form and a homotopy invariant of the
superconnection $\Af$.
\end{prop}
\begin{proof}
Let $\Af_{\s}$ be a 1-parameter family of superconnections on
$\pi_*(\Ee)$ adapted to $\Ps$. Then using the identity
\begin{equation}\label{e:Acommuteswithres}
(\Af_{\s}^2 - \la\Is)\ii \Af_{\s}  = \Af_{\s}(\Af_{\s}^2 -
\la\Is)\ii
\end{equation}
we have formally
\begin{eqnarray}
\dd_{\s} \str(\dd_{\la}^{m-1}(\Af_{\s}^2 - \la\Is)\ii) & =  &
\str(\dd_{\la}^{m-1}\dd_{\s} (\Af_{\s}^2 - \la\Is)\ii) \nonumber \\
& =  & - \str\left(\dd_{\la}^{m-1} \left((\Af_{\s}^2 - \la\Is)\ii
(\dot{\Af}_{\s}\Af_{\s} + \Af_{\s}\dot{\Af}_{\s})\ii
(\Af_{\s}^2 - \la\Is)\ii\right)\right)\nonumber \\
& =  & - \str\left(\left[\Af_{\s}, \dd_{\la}^{m-1}((\Af_{\s}^2 -
\la\Is)\ii \dot{\Af}_{\s} (\Af_{\s}^2 -
\la\Is)\ii\right]\right) \nonumber\\
& = & - d \ \str\left(\dd_{\la}^{m-1}\left((\Af_{\s}^2 -
\la\Is)\ii \dot{\Af}_{\s} (\Af_{\s}^2 - \la\Is)\ii\right)\right)
\label{e:dstr} \\
& = & - d \ \str\left(\dot{\Af}_{\s}\,\dd_{\la}^m\left((\Af_{\s}^2
- \la\Is)\ii \right)\right) \label{e:homotopy} \ ,
\end{eqnarray}
where each step is easily justified rigourously using the kernel
$K(\dd_{\la}^{m-1}(\Af_{\s}^2 - \la\Is)\ii)$, which depends
smoothly on $\s$. The equality \eqref{e:dstr} follows by
essentially the same argument as that in \cite{BGV} Lemma (9.15),
using a parametrix for $\Ps = \Af_{[0]}$ to see the supertrace
vanishes on the supercommutators of vertical $\pdo$s that arise in
\eqref{e:dstr}. (Generally, if $\Rs\in\Aa(B,\Psi(\Ee))$ has
sufficiently negative $\pdo$ order, then one has $d\,\str(\Rs) =
\str([\Af,\Rs])$, but in \eqref{e:dstr} no such subtleties enter.)
Since the variation of the resolvent trace form is exact, this
proves the homotopy invariance.

Similarly, we obtain
\begin{equation}\label{e:dstrclosed}
d \ \str(\dd_{\la}^{m-1}(\Af^2 - \la\Is)\ii) = \str([\Af,
\dd_{\la}^{m-1}(\Af^2 - \la\Is)\ii)]
\end{equation}
which, since the resolvent form has even parity, vanishes by
\eqref{e:Acommuteswithres}, proving closure.
\end{proof}

 From here on we restrict our attention to
the scaled superconnection $\Af_t := t^{1/2}\d_t(\Af)$ with
curvature $\Fs_t := \Af^2_t  = t\d_t(\Af^2) \in \Aa(B,\Psi(\Ee)).
$

\vskip2mm

The small time asymptotics \eqref{e:strrescaled} of the resolvent
trace form for the superconnection curvature and equation
\eqref{e:dstrclosed} now yield the next result, given in terms of
the corresponding coefficient forms in the zeta-form singularity
structure \eqref{e:singularitystructure}.

\begin{prop}\label{p:closedexact}
For
$$ j\neq  r + \frac{rd}{2} + w_k - k + n$$
the $\Ci$ differential forms $a_{j,[d,k]}$ are exact, for $$ l\neq
1 + \frac{d}{2} - q$$ the  forms $a_{l,[d,q]}^{'} ,
a_{l,[d,q]}^{''}$ are exact. The forms
\begin{equation}\label{e:closedexact1}
a_{r + \frac{rd}{2} + w_k - k + n,[d,k]} \ , \hskip 5mm a_{1 +
\frac{d}{2} - q,[d,q]}^{'} \ , \hskip 5mm a_{1 + \frac{d}{2} -
q,[d,q]}^{''} \ ,
\end{equation}
are closed in $\Aa(B)$.

\vskip 2mm Similarly, if $\nu_k\in\Nf$, if
$$ j\neq r + \frac{rd}{2} + w + n$$
the $\Ci$ differential forms $b_{j,d}$ are exact, for $$ l\neq 1 +
\frac{d}{2}$$ the  forms $b^{'}_{l,d}, b^{''}_{l,d}$ are exact.
The forms
\begin{equation}\label{e:closedexact2}
b_{r + \frac{rd}{2} + w + n,d} \ , \hskip 5mm b_{ 1 +
\frac{d}{2},d}^{'} \ , \hskip 5mm b_{ 1 + \frac{d}{2},d}^{''} \ ,
\end{equation}
are closed in $\Aa(B)$.
\end{prop}

For the large time asymptotics we assume that $\Ps = \Fs_{[0]}$
satisfies the constant kernel dimension condition, so that we have
the finite-rank superbundle $\Ker(\Ps)$ over $B$ with the induced
connection $\nabla_0 = \Pi_0\cdot \Af_{[1]}\cdot\Pi_0$ (in the
usual sense), as in \thmref{t:1}, and so we have the corresponding
classical resolvent trace form
\begin{equation*}
    \str\left((\nabla_0^2 - \la I)\ii\right) \in \Aa(B)\ .
\end{equation*}

\begin{prop}\label{p:resformatinfinity}
The following limit holds
\begin{equation}\label{e:resformatinfinity}
\lim_{t\to \o}\str\left(\dd_{\la}^m(\Fs_t - \la\Is)\ii\right) =
\dd_{\la}^m \str\left((\nabla_0^2 - \la I)\ii\right)
\end{equation}
in each $C^l$ norm on compact subsets of $B$. For $t>0$ one has in
$\Aa(B)$
\begin{equation}\label{e:restransgression}
\dd_{\la}^m \str\left((\nabla_0^2 - \la I)\ii\right) =
\str\left(\dd_{\la}^m(\Fs_t - \la\Is)\ii\right) -  d \ \int_t^{\o}
\str\left(\dot{\Af}_{s}\,\dd_{\la}^m\left((\Af_{s}^2 - \la\Is)\ii
\right)\right)\, ds  \ .
\end{equation}
\end{prop}
\begin{proof} We follow the method of \cite{BGV}, Corollary 9.32.
to see that
\begin{eqnarray}
(\Fs_t - \la\Is)\ii & = & \begin{bmatrix}
 (\nabla_0^2 - \la I)\ii  & 0 \\
  0 & 0 \\
\end{bmatrix} + \begin{bmatrix}
 O(t^{-1/2})  & O(t^{-1/2}) \\
  O(t^{-1/2}) & O(t^{-1}) \\
\end{bmatrix} \label{e:testimate}
\end{eqnarray}
as $t\to \o$, and hence that for $m
> (w+n)/r$ in each $C^l$ norm the kernel estimate
$$\|K(\dd_{\la}^m(\Fs_t - \la\Is)\ii) -
K(\dd_{\la}^m (\nabla_0^2 - \la I)\ii) \|_l \leq C_l t^{-1/2}$$
holds uniformly on compact subsets of $M\times_{\pi}M$, from which
\eqref{e:resformatinfinity} follows. By integrating
\eqref{e:homotopy} we have for $0 < t < T < \o$
\begin{equation}\label{e:transgression1}
\str(\dd_{\la}^{m-1}(\Fs_T - \la\Is)\ii) -
\str(\dd_{\la}^{m-1}(\Fs_t - \la\Is)\ii) = - d \ \int_t^T
\str\left(\dot{\Af}_{s}\,\dd_{\la}^m\left((\Af_{s}^2 - \la\Is)\ii
\right)\right)\, ds  \ .
\end{equation}
On the other hand, the estimate \eqref{e:testimate} implies that
\begin{equation}\label{e:strestimate2}
\|\str(\dot{\Af}_{s}\,\dd_{\la}^m ((\Af_{s}^2 - \la\Is)\ii )) \|_l
= O(t^{-3/2})
\end{equation}
as $t\to\o$ in each $C^l$ norm on compact subsets of $B$, and so
with \eqref{e:resformatinfinity} the identity
\eqref{e:restransgression} follows.
\end{proof}

\begin{cor}\label{c:restransgression}
The following cohomology identity holds. For $m
> (w+n)/r$ in $H^*(B)$
\begin{equation}\label{e:cohomologyrestransgression}
\dd_{\la}^m \str\left((\nabla_0^2 - \la I)\ii\right) = \LIM_{t\to
0}\str\left(\dd_{\la}^m(\Fs_t - \la\Is)\ii\right)
\end{equation}
\begin{equation}\label{e:cohomologyrestransgression2}
= \sum_{d,k=0}^{\dim B}(\a_{w_k + n -r k + r + \frac{rd}{2},[d,k]}
+ \a^{'}_{1-k + \frac{d}{2},[d,k]}\log (\la) + \a^{''}_{1-k +
\frac{d}{2},[d,k]})\, (-\la)^{-1-\frac{d}{2}-m} \ ,
\end{equation}
where \eqref{e:cohomologyrestransgression2} follows from
\eqref{e:strrescaled}, and, at the level of differential forms,
the coefficients are closed, differing from the forms in
\eqref{e:closedexact1} by constants.
\end{cor}

\vskip 3mm

Notice that \eqref{e:heattrace2} and \eqref{e:restransgression}
prove the Chern character transgression formula of \cite{B,BGV}.
Indeed, the above formulas are the governing transgression
formulas for all characteristic class forms on $\pi_*(\Ee)$.

\section{Zeta Forms and the Family Index Theorem}

This Section consists of the proof of \thmref{t:1}. Throughout
\begin{equation}\label{e:rescaledconnectionandcurvature}
\Af_t := t^{1/2}\d_t(\Af) \ , \hskip 10mm \Fs_t := \Af^2_t  =
t\d_t(\Af^2) \ .
\end{equation}
Evidently $\Af =  \Af_1, \Fs = \Fs_1$.

\vskip 3mm
\begin{prop}\label{p:zetaformexact}
The zeta form $\zeta(\Af^2,s)|^{{\rm mer}}$ is canonically exact.
One has in $\Aa(B)$
\begin{equation}\label{e:zetaformexact}
    \zeta(\Af^2,s)|^{{\rm mer}} = d\ \int^{\o}_1
    \zeta(\dot{\Af}_{\s}. \Fs_{\s} , s) |^{{\rm mer}} \, d\s \ .
\end{equation}
\end{prop}
\begin{proof}
For simplicity we assume \eqref{e:orderlessthanr} holds, the
modifications for the general case are obvious.

For $\re(s)>0$ we then have
\begin{equation*}
\Fs_{\th}\si = \frac{i}{2\pi}\int_{C}\la_{\th}\si\, (\Fs -
\la\Is)\ii \, d\la \ .
\end{equation*}
It follows that on $\Ker(\Ps)$
\begin{equation*}
\Fs\si_{|{\tiny \Ker(\Ps)}} \equiv 0 \ .
\end{equation*}
For, from \eqref{e:resexpansion0}
\begin{equation}\label{e:resonker}
(\Fs - \la \Is)\ii _{|{\tiny \Ker(\Ps)}} = - \sum_{i=0}^{\dim B}
\la^{-i-1}\left(\Pi_0\cdot\Ws\cdot\Pi_0\right)^i \ ,
\end{equation}
and $\frac{i}{2\pi}\int_{C}\la_{\th}^{-s-i-1} \, d\la = 0$ for
$i\geq 0$ and $\re(s) >0$. Hence we have
\begin{equation}\label{e:vanishesonkernel}
\str(\Fs\si_{|{\tiny \Ker}(\Ps)}) = 0 \ , \hskip 10mm \re(s) > 0 \
.
\end{equation}
On the other hand, \eqref{e:testimate} gives
\begin{equation*}
\str\left((\Fs_T - \la \Is)\ii _{|{\tiny \Ker(\Ps)}\pp}\right) =
O(T\ii)\hskip 5mm {\rm as} \ \ T\too \o \ .
\end{equation*}
Consequently from \eqref{e:transgression1}
\begin{equation}\label{e:strexact}
\str(\dd_{\la}^{m-1}(\Fs - \la\Is)\ii_{|{\tiny \Ker(\Ps)}\pp}) =
d \ \int_1^{\o}
\str\left(\dot{\Af}_{\s}\,\dd_{\la}^m\left((\Af_{\s}^2 -
\la\Is)\ii \right)\right)\, d\s  \ .
\end{equation}
From \eqref{e:vanishesonkernel}, \eqref{e:strexact} and
\eqref{e:superzetaform2} we find for $\re(s) > m > (w+n)/r$
\begin{equation*}%\label{e:slargezetaexact}
\z_{\th}(\Fs,s)
\end{equation*}
$$= d \ \left( \int_1^{\o} \frac{1}{(s-1)\ldots
(s-m)}.\frac{i}{2\pi}\int_{C}\la^{m-s}\,
\str\left(\dot{\Af}_{\s}\,\dd_{\la}^m (\Af_{\s}^2 - \la\Is)\ii
\right) \, d\la  \, d\s\right)$$
$$= d\ \int^{\o}_1
    \zeta(\dot{\Af}_{\s}. \Fs_{\s} , s) \, d\s \ . $$
Hence in $\Aa(B)$
\begin{equation}\label{e:slargezetaexact}
\z_{\th}(\Fs,s) - d\ \int^{\o}_1
    \zeta(\dot{\Af}_{\s}. \Fs_{\s} , s) \, d\s  = 0 \ ,\hskip 10mm
    \re(s) > (w+n)/r \ .
\end{equation}

Elsewhere in $\Cf$, from \eqref{e:resonker} and \cite{GS1}
Proposition (2.9), we see that $\G(s)\str(\Fs\si_{|{\tiny
\Ker}(\Ps)})$ has no poles. (Less strong, but more general, and
sufficient, the zeta form for any family of finite rank operators
extends without poles.) Hence $\str(\Fs\si_{|{\tiny
\Ker}(\Ps)})|^{{\rm mer}} $ is a holomorphic extension of zero to
all of $\Cf$, and consequently
\begin{equation}\label{e:vanishesonkernelforalls}
\str(\Fs\si_{|{\tiny \Ker}(\Ps)})|^{{\rm mer}} = 0  \ .
\end{equation}
Likewise
\begin{equation}\label{e:exactformforalls}
(\, d\,\int^{\o}_1 \zeta(\dot{\Af}_{\s}. \Fs_{\s} , s)\,d\s \
)|^{{\rm mer}} = d\ (\int^{\o}_1 \zeta(\dot{\Af}_{\s}. \Fs_{\s} ,
s)|^{{\rm mer}} \, d\s \ )  \ ,
\end{equation}
since from \eqref{e:homotopy} and \eqref{e:singularitystructure}
both sides of \eqref{e:exactformforalls} have the same pole
structure. Hence we find that $\z_{\th}(\Fs,s)|^{{\rm mer}} - d\
\int^{\o}_1 \zeta(\dot{\Af}_{\s}. \Fs_{\s} , s)|^{{\rm mer}} \,
d\s $ is holomorphic on $\Cf$ and so from
\eqref{e:slargezetaexact} it is identically zero.
\end{proof}

Evidently, then, \eqref{e:familyzetaindex} now follows with
\begin{equation}\label{e:familyzetaindex3}
\sum_{k=0}^{\dim B}\frac{1}{k!}\ \z_{\pi}(\Fs,-k)|^{{\rm mer}} = d
\left(\sum_{k=0}^{\dim B}\frac{1}{k!}\ \int^{\o}_1
\zeta(\dot{\Af}_{\s}. \Fs_{\s} , -k)|^{{\rm mer}} \, d\s \ \right)
\ .
\end{equation}

\vskip 2mm

For clarity, and as it applies to the case of the Bismut
connection which is the superconnection of primary geometric
interest, we restrict our formulas from here on to the case
\begin{equation}\label{e:integer order}
    \nu_i = \ord(\Fs_{[i]}) \in \Nf \ .
\end{equation}
The case for any real $\nu_i$ is the same but with more indices to
track.

\vskip 3mm

From \eqref{e:strexpansion2} the $t$-rescaled singularity
structure equation \eqref{e:singularitystructure2} is
\begin{equation}\label{e:singularitystructure2rescaled}
 \frac{\pi}{\sin(\pi s)}\z_{\th}(\Fs_t,s)|^{{\rm mer}} \sim
 \end{equation}
 $$-\sum_{j=-\dim B -1}^{-1}
  \frac{\str\left((\Pi\cdot\Ws_t\cdot\Pi)^{-j-1}\right)}{(s - j -
  1)^{k+1}} \ \ + $$
 $$\sum_{d=0}^{\dim B}\left(\sum_{j\geq 0}
\frac{b_{j,d}}{\left(s + \frac{j - n - w}{r} -1 \right)}
t^{\frac{j - n - w}{r} -1 - \frac{d}{2}}\right. + \left.
\sum_{l\geq 0} \frac{b^{''}_{l,d}}{(s + l - 1)}t^{l -1 -
\frac{d}{2}}\right) \ ,$$ where
\begin{equation}\label{e:Wt}
    \Ws_t = t\d_t(\Ws) = \Fs_t - t\Ps^2 \ .
\end{equation}

For $s$ in small neighborhood of $-k$ there is a Laurent expansion
\begin{equation*}
\frac{\sin(\pi s)}{\pi} = (-1)^k (s-k) + O((s-k)^3) \ ,
\end{equation*}
and hence from  \eqref{e:singularitystructure2rescaled}
\begin{equation}\label{e:zetaat-k}
(-1)^k  \z_{\pi}(\Fs_t,-k)|^{{\rm mer}} =
 - \str\left((\Pi\cdot\Ws\cdot\Pi)^k\right)
 + \sum_{d=0}^{\dim B} \left(b_{n+ w + r + rk,d} +
 b^{''}_{k+1,d}\right)t^{k - \frac{d}{2}} \ .
\end{equation}
On the other hand, \eqref{e:closedexact2} of
\propref{p:closedexact} says that the forms $b_{n+ w + r + rk,d}$
are exact except possibly when
$$n+ w + r + r k = r + \frac{rd}{2} + w + n \ .$$
That is, when $d = 2k.$ Like wise for the $b^{''}_{k+1,d}$, and so
\eqref{e:zetaat-k} can be written \vskip 2mm
\begin{equation*}\label{e:zetaat-k2}
\z_{\pi}(\Fs_t,-k)|^{{\rm mer}} =
-(-1)^k\str\left((\Pi\cdot\Ws_t\cdot\Pi)^k\right) +  b_{n+ w + r +
rk,2k} + b^{''}_{k+1,2k} + d \gamma^{t}_{n,w,r,k}\ ,
\end{equation*}
\vskip 3mm \noi where $\gamma^t_{n,w,r,k}\in\Aa(B)$ such that
$d\gamma^t_{n,w,r,k}$ is the sum of exact forms in
\eqref{e:zetaat-k} minus the two closed forms above. Hence
\begin{eqnarray}
\sum_{k=0}^{\dim B}\frac{1}{k!}\ \z_{\pi}(\Fs_t,-k)|^{{\rm mer}} =
& - & \sum_{k=0}^{\dim B}\frac{(-1)}{k!}^k
\str\left((\Pi\cdot\Ws_t\cdot\Pi)^k\right) \nonumber \\
& & + \sum_{k=0}^{\dim B}\frac{1}{k!}\,(b_{n+ w + r + r k,2k} +
b^{''}_{k+1,2k}) + d \sum_{k=0}^{\dim
B}\frac{1}{k!}\,\gamma^t_{n,w,r,k}  \nonumber
\end{eqnarray}
\hskip 25mm ---  which since
$\str\left((\Pi\cdot\Ws_t\cdot\Pi)^k\right)
 \in\Aa^{>0}(B)$ for $k>0$ ---
\begin{eqnarray}
\hskip 39mm = & - & \str\left(e^{-(\Pi\cdot\Ws_t\cdot\Pi)^k}\right) \nonumber\\
& + & \sum_{k=0}^{\dim B}\frac{1}{k!}\,(b_{n+ w + r + r k,2k} +
b^{''}_{k+1,2k}) + d \sum_{k=0}^{\dim
B}\frac{1}{k!}\,\gamma^t_{n,w,r,k} \  . \label{e:zetasum}
\end{eqnarray}

\noi Now from \eqref{e:zetaat-k} we have $\LIM_{t\to 0} \ d
\gamma^{t}_{n,w,r,k} = 0.$ On the other hand, the $t$-independent
component of $\Ws_t$ is the 2-form piece which is the curvature
form $\Af_{[1]}^2$. Hence we find $\LIM_{t\to 0} \
\str\left((\Pi\cdot\Ws_t\cdot\Pi)^k\right)  =
\str\left(\nabla_0^{2k}\right),$ with $\nabla_0
=\Pi\cdot\Af_{[1]}\cdot\Pi$ the induced connection on the
superbundle $\Ker(\Ps)$. Consequently, since
\eqref{e:heatraceexpansion2} and \eqref{e:heatzetafactors} give
\begin{equation*}\
\LIM_{t\to 0} \ \str(e^{-\Fs_t}) = \sum_{k=0}^{\dim
B}\frac{1}{k!}\,(b_{n+ w + r + r k,2k} + b^{''}_{k+1,2k}) \ ,
\end{equation*}
taking the regularized limit as $t\to 0+$ of \eqref{e:zetasum} we
find
\begin{equation}\label{e:chcharacter1}
\LIM_{t\to 0} \ \sum_{k=0}^{\dim B}\frac{1}{k!}\
\z_{\pi}(\Fs_t,-k)|^{{\rm mer}} = - \ch(\Ker(\Ps),\nabla_0) +
\LIM_{t\to 0} \ \ch(\Af_t) \ .
\end{equation}

\noi Since, from \eqref{e:familyzetaindex3}, the left-side of
\eqref{e:chcharacter1} is equal to
\begin{equation}\label{e:familyzetaindex4}
d \ \LIM_{t\to 0} \int^{\o}_t \str(\, \dot{\Af}_{\s}
\sum_{k=0}^{\dim B} \frac{(-1)}{k!}^k \Fs_{\s}\si|_{s=-k}^{{\rm
mer}}\ )\ \ d\s
\end{equation}
$$ = d \ \LIM_{t\to 0}  \int^{\o}_t \str(\, \dot{\Af}_{\s}e^{-\Fs_{\s}}
 )\ \ d\s \ $$
\vskip 2mm \noi this completes the proof of \thmref{t:1}.

\section{Zeta-Chern Forms}

In this Section we prove \thmref{t:2}.

\begin{defn}
Let $\Af$ be a superconnection on $\p_*(\Ee)$ adapted to a family
of formally self-adjoint  $\z$-admissible $\pdo$s
$\Ps\in\G(B,\Psi^{r>0}(\Ee))$ of odd-parity. The zeta Chern form
of $\Af$ is the differential form of mixed order on $B$ defined by
the super zeta determinant form
\begin{equation}\label{e:zetachernform}
    \cc_{\, \z}(\Af) = {\rm sdet}_{\z,\pi}(\Is + \Af^2) \in
    \Aa(B)\ .
\end{equation}
Here $\Is\in\Aa^0(B,\Psi^{0}(\Ee))$ is the vertical identity
operator and $\Af^2\in\Aa(B,\Psi(\Ee))$ the superconnection
curvature.
\end{defn}

\vskip 1mm

Thus, with $ \Fs = \Af^2$ we have
\begin{equation}\label{e:logzetachernform}
-\log\cc_{\, \z}(\Af) = \dd_s \str\left((\Is +
\Fs^2)\si\right)|^{{\rm mer}}_{s=0} = \str\left(\log(\Is +
\Fs)\,(\Is + \Fs)\si \right)|^{{\rm mer}}_{s=0} \ .
\end{equation}

\vskip 2mm Notice  that since
\begin{equation}\label{e:0formpart}
(\Is + \Fs)_{[0]} = \Is + \Ps^2
\end{equation}
then $\Is + \Fs\in \Aa(B,\Psi^{r>0}(\Ee))$ is elliptic and
invertible with Agmon angle $\th = \pi$. Further, since
\begin{equation*}
{\rm sdet}_{\z,\pi}(\Is + \Ps^2) = 1 \ ,
\end{equation*}
then by \lemref{p:zetaperturbation} we have $\log\cc_{\,
\z}(\Af)_{[0]} = 0$ and hence
\begin{equation}\label{e:0formpartofzetachern}
\cc_{\, \z}(\Af)_{[0]} = 1 \ .
\end{equation}

To write down the singularity structure of the zeta-form
$\z_{\pi}(\Is + \Fs,s)$ it is convenient to use the function
introduced in \cite{GS2} defined for $\re(-t) < \re(s) <0$ by
\begin{equation}\label{e:Fts}
    F_t(s) = \frac{i}{2\pi}\int_{C} \mu^{-s-1}(1-\mu)^{-t} \, d\mu
\end{equation}
with $C$ a contour around $R_{\pi}$. $F_t(s)$ extends
meromorphically to all of $\Cf$ and satisfies
\begin{equation}\label{e:Ftsidentities}
    F_t(s)  =  \frac{\G(s+t)}{\G(t)\G(s+1)} .
\end{equation}

Since
\begin{equation*}
\str\left(\dd_{\la}^{m-1}((\Is + \Fs) - \la \Is)\ii\right) =
\str\left(\dd_{\la}^{m-1}(\Fs + (1-\la)\Is)\ii\right)
\end{equation*}
we have from \eqref{e:strexpansion2} an asymptotic expansion of
the resolvent trace as $\la\to\o$ in $\Lambda_{\pi}$
\begin{eqnarray}
\str(\partial^{m-1}_{\la}(\Fs - (1-\la)\ \Is)\ii\,)  & = &
\sum_{d=0}^{\dim B}\left(\sum_{j=0}^{N-1} (-1)^{\frac{w + n
-j}{r}}\b_{j,d}(1-\la)^{\frac{w + n -j}{r}-m}\right. \nonumber \\
& & \left. + \sum_{l=0}^{N-1}  (-1)^{-l}(\b^{'}_{l,d}\log (1-\la)
+ \b^{''}_{l,d})(1-\la)^{-l-m}\right) \nonumber \\
& & + \ O(| 1 -\la| ^{^{\frac{w + n -N}{r}-m}}) \ .
\label{e:resstrdet}
\end{eqnarray}
Hence in $\Aa(B)$
\begin{eqnarray}
\str((\Is + \Fs)\si\,)|^{{\rm mer}}  & = & \sum_{d=0}^{\dim
B}\left(\sum_{j=0}^{N-1} (-1)^{\frac{w + n
-j}{r}} b_{j,d} \ F_{\frac{j-w -n}{r}}(s-1)|^{{\rm mer}}\right. \nonumber \\
& & \left. + \sum_{l=0}^{N-1} (-1)^{-l} b^{'}_{l,d} \ \dd_s
F_l(s-1)|^{{\rm mer}}
+ \sum_{l=0}^{N-1} b^{''}_{l,d} \ F_l(s-1)|^{{\rm mer}}\right)  \nonumber \\
& & + \ h_N(s) \ , \label{e:singstrdet}
\end{eqnarray}
where $h_N(s)\in \Aa(B)$ is holomorphic for
\begin{equation}\label{e:hol}
    1 - \left(\frac{N-n-w}{r}\right) < \re(s) < N+1 \ ,
\end{equation}
and $b_{j,d}, b^{'}_{l,[d,q]}, b^{''}_{l,[d,q]}\in\Aa^{d}(B)$.

\vskip 5mm

\begin{prop}\label{p:zetachernclosed}
The differential form $\cc_{\, \z}(\Af)$ is closed.
\end{prop}
\begin{proof}
For $\re(s) >> 0$ we have
\begin{equation*}%\label{e:detreslarge}
    \str\left(\log(\Is +
\Fs)\,(\Is + \Fs)\si \right) = \frac{1}{(s-m)\ldots(s-1)}
\frac{i}{2\pi}\int_C \log \la \, \la^{m-s} \,
\str(\dd_{\la}^{m-1}(\Fs + (1-\la)\Is)\ii) \ d\la
\end{equation*}
and hence from \propref{p:resformclosed}
\begin{equation}\label{e:ddetzerolarges}
    d \ \str\left(\log(\Is +
\Fs)\,(\Is + \Fs)\si \right) = 0 \ , \hskip 15mm \re(s) >> 0 \ .
\end{equation}
Elsewhere from \eqref{e:singstrdet}
\begin{eqnarray}
\str\left(\log(\Is + \Fs)\,(\Is + \Fs)\si \right)|^{{\rm mer}}  &
= & \sum_{d=0}^{\dim B}\left(\sum_{j=0}^{N-1} (-1)^{\frac{w + n
-j}{r}} b_{j,d} \ \dd_s F_{\frac{j-w -n}{r}}(s-1)|^{{\rm mer}}\right. \nonumber \\
& & \left. + \sum_{l=0}^{N-1} (-1)^{-l} b^{'}_{l,d} \ \dd^2_s
F_l(s-1)|^{{\rm mer}}
+ \sum_{l=0}^{N-1} b^{''}_{l,d} \ \ \dd_s F_l(s-1)|^{{\rm mer}}\right)  \nonumber \\
& & + \ \dd_s h_N(s) \ ,  \label{e:singstrdet1}
\end{eqnarray}
But from \propref{p:closedexact} the forms $b_{j,d},
b^{'}_{l,[d,q]}, b^{''}_{l,[d,q]}$ are all closed. Hence
\begin{equation}\label{e:ddetreslarge}
  d \   \str\left(\log(\Is + \Fs)\,(\Is + \Fs)\si \right)|^{{\rm mer}}
 = d  \ \dd_s h_N(s) \ .
\end{equation}
The right-side of \eqref{e:ddetreslarge} is independent of $N$ and
holomorphic. By \eqref{e:ddetzerolarges} we obtain that $d \
\str\left(\log(\Is + \Fs)\,(\Is + \Fs)\si \right)|^{{\rm mer}}$ is
a holomorphic extension of zero and hence vanishes identically on
all of $\Cf$, proving the assertion.
\end{proof}

The zeta-Chern form $\cc_{\, \z}(\Af)$ hence defines a mixed
degree cohomology class in $H^*(B)$. To see this is the Chern
class of the index bundle we have the following transgression
results.

\vskip 5mm

\begin{prop}\label{p:zetacherntransgression1}
If $\Af_{\s}$ is a 1-parameter family of superconnections adapted
to $\Ps\in\Aa^0(B,\Psi(\Ee))$ with curvature $\Fs_{\s} =
\Af^2_{\s}$, then
\begin{equation}\label{e:zetacherntransgression1}
\dd_{\s} \log \cc_{\, \z}(\Af_{\s}) = - d \ \zeta( (\Fs_{\s} +
\Is)\ii \dot{\Af}_{\s}, \Fs_{\s} + \Is , 0)|^{{\rm mer}} \ .
\end{equation}
Here, $\zeta( (\Fs_{\s} + \Is)\ii \dot{\Af}_{\s}, \Fs_{\s} + \Is ,
s)|^{{\rm mer}} = \str((\Fs_{\s} + \Is)^{-s-1}
\dot{\Af}_{\s})|^{{\rm mer}}$.
\end{prop}
\begin{proof}
From \eqref{e:homotopy} we have for sufficiently large $m$
\begin{equation}\label{e:restrvariation}
\dd_{\s} \ \str(\dd_{\la}^{m-1}(\Fs_{\s} + (1-\la)\Is)\ii) = -d \
\str(\dd_{\la}^m(\Fs_{\s} + (1-\la)\Is)\ii \dot{\Af}_{\s} ) \ .
\end{equation}
Hence for $\re(s) >> 0$ we have integrating by parts
\begin{equation}\label{e:restrvariation2}
\dd_{\s}\ \str\left(\log(\Is + \Fs_{\s})\,(\Is + \Fs_{\s})\si
\right)
\end{equation}
\begin{eqnarray*}
& = & -d \ \left(\frac{1}{(s-m)\ldots(s-1)} \frac{i}{2\pi}\int_C
\log \la \, \la^{m-s} \, \str(\dd_{\la}^m(\Fs_{\s} +
(1-\la)\Is)\ii
\dot{\Af}_{\s}) \ d\la \right) \nonumber \\
& = & -d \ \left(\frac{1}{(s-m)\ldots(s-1)} \frac{i}{2\pi}\int_C
\la^{m-s-1} \, \str(\dd_{\la}^{m-1}(\Fs_{\s} + (1-\la)\Is)\ii
\dot{\Af}_{\s}) \ d\la \right) \nonumber  \\
& &  + \  s\,d \ \left(\frac{1}{(s-m)\ldots(s-1)}
\frac{i}{2\pi}\int_C  \log\la\,\la^{m-s-1} \,
\str(\dd_{\la}^{m-1}(\Fs_{\s} + (1-\la)\Is)\ii
\dot{\Af}_{\s}) \ d\la \right) \nonumber  \\
& = &- d \ \zeta( (\Fs_{\s} + \Is)\ii \dot{\Af}_{\s}, \Fs_{\s}+
\Is, s) + s \,d \ \dd_s \zeta((\Fs_{\s} + \Is)\ii \dot{\Af}_{\s},
\Fs_{\s}) + \Is , s) \ .
\end{eqnarray*}
Since $\zeta((\Fs_{\s} + \Is)\ii \dot{\Af}_{\s}, \Fs_{\s}) + \Is ,
s)$ is holomorphic around $s=0$, as we consider $\z$-admissible
vertical operators, then evaluating at zero
\eqref{e:restrvariation2} gives
\begin{equation*}%\label{e:restrvariation3}
\dd_{\s} \  \str\left(\log(\Is + \Fs_{\s})\,(\Is + \Fs_{\s})\si
\right)|^{{\rm mer}}_{s=0} = -\left(d \ \zeta( (\Fs_{\s} + \Is)\ii
\dot{\Af}_{\s}, \Fs_{\s}+ \Is, s)\right)|^{{\rm mer}}_{s=0} \ .
\end{equation*}
Since the singularity expansion \eqref{e:singstrdet1} shows that
the derivatives can be commuted with $|^{{\rm mer}}$, this
completes the proof.
\end{proof}

We now restrict attention to the rescaled superconnection and
curvature \eqref{e:rescaledconnectionandcurvature}.

\vskip 5mm

\begin{prop}\label{p:zetacherntransgression2}
Let  $\tau_{t,T}(\Af) = - d \ \int_t^T \zeta( (\Fs_{\e} + \Is)\ii
\dot{\Af}_{\e}, \Fs_{\e} + \Is , 0)|^{{\rm mer}} d\e$ with $0 < t
< T < +\o$. Then in $\Aa(B)$ one has
\begin{equation}\label{e:zetacherntransgression2}
\frac{\cc_{\, \z}(\Af_T)}{\cc_{\, \z}(\Af_t)} = e^{\tau_{t,T}(\Af)
} \ .
\end{equation}
Equivalently,
\begin{equation}\label{e:zetacherntransgression4}
\cc_{\, \z}(\Af_T) = \cc_{\, \z}(\Af_t) + d\, \omega_{t,T}(\Af) \
,
\end{equation}
where
\begin{equation}\label{e:zetacherntransgression6}
\omega_{t,T}(\Af) = \cc_{\, \z}(\Af_t) \wedge \sum_{k\geq 1}
\frac{1}{k!} \tau_{t,T}(\Af) \wedge ( d\,\tau_{t,T}(\Af))^{k-1} \
.
\end{equation}
If $\Ps\in\Aa^0(B,\Psi(\Ee))$ has constant kernel dimension, then
the limit $\lim_{T\to\o}\omega_{t,T}(\Af) $, denoted
$\omega_{t,\o}(\Af)$ exists in all $C^l$-norms on compact subsets
of $B$, and one has
\begin{equation}\label{e:zetacherntransgression5}
\cc(\Ker(\Ps), \nabla_0) = \cc_{\, \z}(\Af_t) + d\
\omega_{t,\o}(\Af) \ .
\end{equation}
with notation as in \thmref{t:2}.
\end{prop}
\begin{proof}
From \eqref{e:0formpartofzetachern} the quotient on the left-side
of \eqref{e:zetacherntransgression2} is well-defined in $\Aa(B)$.
The identity is immediate from integrating
\eqref{e:zetacherntransgression1} and exponentiating both sides.
Since $\cc_{\, \z}(\Af_t)$ is closed,
\eqref{e:zetacherntransgression4} is immediate.

By \propref{p:resformatinfinity} and \eqref{e:singstrdet1} we find
\begin{eqnarray*}
\lim_{T\to\o} \cc_{\, \z}(\Af_T) & = & \frac{1}{(s-m)\ldots(s-1)}
\frac{i}{2\pi}\int_C \log \la \, \la^{m-s} \,
\lim_{T\to\o}\str(\dd_{\la}^{m-1} (\Fs_T + (1-\la)\Is)\ii  \ d\la
|^{{\rm
mer}}_{s=0} \nonumber  \\
& = & \frac{i}{2\pi}\int_C \log \la \, \str((\nabla_0^2 +
(1-\la)\Is)\ii \ d\la  \nonumber  =  \log \ {\rm sdet}(I +
\nabla_0^2)) \ .
\end{eqnarray*}
Since  \eqref{e:strestimate2}implies the $\z$-form $C^l$ estimate
$ \| \zeta( (\Fs_{\e} + \Is)\ii \dot{\Af}_{\e}, \Fs_{\e} + \Is ,
0)|^{{\rm mer}} \|_l \leq c(l) t^{-3/2}$ we obtain the existence
of the limit $\lim_{T\to\o}\tau_{t,T}(\Af)$ and hence the limit
$\lim_{T\to\o}\omega_{t,T}(\Af)$.
\end{proof}

We turn next to the proof of the local index density formula
\eqref{e:chernlimit}.

\vskip 3mm

We suppose that the fibre bundle $\pi : M \too B$ has
even-dimensional fibre and that it is endowed with a connection
\begin{equation}\label{e:bconnection}
TM = \pi^*(TB)\oplus T(M/B) \ ,
\end{equation}
defined by a choice of bundle projection $P : TM \to T(M/B)$.
Suppose also that $TM$ has a spin structure and that there are
Riemannian metrics $g_{M/B}, g_B$ on $T(M/B)$ and $TB$. The
vertical bundle $\Ee$ is assumed to be a bundle of Clifford
modules equipped with a connection which restricts to a Clifford
connection on $\Ee_{|M_z}$. Let $\Ds$ be the associated family of
compatible Dirac operators, let $\nabla^{\pi_*(\Ee)}$ be the
canonical Hermitian connection induced on $\pi_*(\Ee)$ \cite{B},
\cite{BGV} Proposition(9.13), and let $\cc(T)\in
\Aa^2(B,\End(\pi_*(\Ee)))$  denote Clifford multiplication by the
torsion tensor of the fibration associated to the connection
\eqref{e:bconnection}. For $t>0$ the scaled Bismut superconnection
on $\pi_*(\Ee)$ is then defined by \cite{B,BGV}
\begin{equation}\label{e:bsuperconnection}
\Af_t = t^{1/2}\Ds + \nabla^{\pi_*(\Ee)} + \frac{1}{4t^{1/2}}\,
\cc(T) \ .
\end{equation}
Its crucial property \cite{BGV} Proposition (10.28) is that with
$\Fs_t = \Af_t^2\in \Aa^4(B,\Psi(\Ee))$, in a small enough
neighborhood $U$ of $x_0\in M$ as the origin of local geodesic
coordinates $x = (x_1,\ldots,x_n)$ along the fibres, the limit
$\lim_{t\to 0}\str(e^{-\Fs_t})$ exists and is equal to the heat
kernel trace $\str(e^{-L_{\xx}})|_{\textbf{x}=0}$, where $\xx\in
T_xM$ has norm $\|\xx\|_{g}$ less than the injectivity radius of
$M$, so that $x = \exp_{x_0}\xx$, of the localized operator
\begin{equation*}\label{e:smalltbsuperconnection}
L_{\xx} = -\sum_i \left(\dd_i -
\frac{1}{4}\sum_j(R^{M/B}\dd_i,\dd_j)\xx_j \right)^2 + R^{\Ee/S} \
.
\end{equation*}
Equivalently
\begin{equation}\label{e:smalltbsuperconnection}
\lim_{t\too 0} \str\left( \dd_{\la}^{m-1}(\Fs_t -
\la\Is)\ii\right) = \str(\partial^{m-1}_{\la}(L_{\xx} - \la
\Is)_{|U}\ii\,)|_{\xx=0}.
\end{equation}
This is demonstrated via Getzler rescaling of the time, fibre and
vertical Clifford variables. Let
$$\widehat{A}(M/B) = {\rm det}^{1/2}\left(\frac{R^{M/B}/2}{\sinh
(R^{M/B}/2)}\right)$$ be the vertical $\widehat{A}$-genus form for
the connection $\nabla^{M/B} = P\cdot \nabla^M \cdot P$, with
$\nabla^M $ the Levi-Civita connection defined by the metric
$g_{M/B}\oplus\pi^*(g_B)$, and let $\ch^{'}(\Ee) =
\str_{\Ee/S}(e^{-R^{\Ee/S}})$ be the relative Chern character of
$\Ee$ and the spin bundle $S$ on $M$.

\begin{prop}\label{p:smalltlimitresstrace}
With $\Fs_t$ the Bismut superconnection curvature, the
differential form has a limit as $t\to 0$ given by the formula
\begin{equation}\label{e:smalltlimitresstrace}
    \lim_{t\too 0} \str\left( \dd_{\la}^{m-1}(\Fs_t -
    \la\Is)\ii\right)
\end{equation}
$$ = (m-1)! \sum_{k=0}^{[\dim B/2]}(-\la)^{-1-k-m}
    \, k! \left((2\pi )^{-\frac{n}{2}}\int_{M/B} \widehat{A}(M/B) \ch^{'}(\Ee)
    \right)_{[2k]} \ .$$
\end{prop}
\begin{proof}
Since the component terms $(\Fs_t)_{[i]}\in\Aa^i(B,\Psi(\Ee))$,
$i=0,\ldots,4,$  are smooth families of {\em differential}
operators, with $(\Fs_t)_{[0]} = t\Ds^2$ so that $r=2$, it follows
from \cite{GS1} that there are no log terms in the resolvent
supertrace and that all the coefficients are local, determined by
only finitely many terms of the vertical symbol.; that is,
$\b^{'}_{l,d} = 0, \b^{''}_{l,d} = 0$ in \eqref{e:strrescaled2},
so that as $t\to 0+$ there is an asymptotic expansion
\begin{equation}\label{e:strrescaled3}
\str(\partial^{m-1}_{\la}(\Fs_t - \la \Is)\ii\,)  \sim
\sum_{d=0}^{\dim B}\sum_{j\geq 0} \b_{j,d}(-\la)^{\frac{w + n
-j}{2}-m} \ t^{\frac{j - w - n}{2} -1 -\frac{d}{2}} \ .
\end{equation}
On the other hand, it follows from
\eqref{e:smalltbsuperconnection} that
\begin{equation}\label{e:strrescaled4}
 \lim_{t\too 0} \str\left( \dd_{\la}^{m-1}(\Fs_t -
    \la\Is)\ii\right)
\end{equation}
$$= \sum_{k=0}^{\dim B}\partial^{m-1}_{\la}\left(\
(\sum_i\dd_i^2-\la\Is)\ii\left(\omega(R^{M/B},R^{\Ee/S})(\xx)\
(\sum_i\dd_i^2-\la\Is)\ii \ \right)^k\right)|_{\xx=0} \ ,$$ with
$\omega(R^{M/B},R^{\Ee/S})(\xx)\in\Aa^{k\geq 2}(B)$ a form of
mixed degree 2 or greater which is $O(1)$ in $\xx$, and hence that
the resolvent supertrace has a limit as $t\to 0$.

Consequently, the expansion \eqref{e:strrescaled3} begins with the
$t^0$ term, that is when
\begin{equation}\label{e:t0}
\frac{j-n-w}{2} = 1 + \frac{d}{2} \ ,
\end{equation}
and since from \eqref{e:strrescaled4} the form degree is always
even  $d=2k$, then as $t\to 0+$
\begin{equation}\label{e:strrescaled6}
\str(\partial^{m-1}_{\la}(\Fs_t - \la \Is)\ii\,)  =
\sum_{k=0}^{[\dim B/2]} \b_{n+w+2+2k,2k}[m](-\la)^{-1-k-m}  +
O(t^{1/2}) \ .
\end{equation}

We hence obtain from \eqref{e:singularitystructure2rescaled}
\begin{equation}\label{e:singularitystructure2rescaled2}
 \frac{\pi}{\sin(\pi s)}\z_{\th}(\Fs_t,s)|^{{\rm mer}} =
 \end{equation}
 $$-\sum_{l=0}^{\dim B}
  \frac{\str\left((\Pi\cdot\Ws_t\cdot\Pi)^l\right)}{(s + l)} \ \ +
 \sum_{k=0}^{[\dim B/2]}\frac{b_{n+w+2+2k,2k}}{\left(s + k\right)} +
 G_{t^{1/2}}(s) \ ,$$ where $G_{t^{1/2}}(s)  = O(t^{1/2})$ and is
 meromorphic on $\Cf$ with no poles at the negative integers, and
 according to \eqref{e:singcoeffs2} and \eqref{e:heatzetafactors}
\begin{equation}\label{e:singcoeffs3}
\b_{n+w+2+2k,2k}[m]  = (m-1)!\ b_{n+w+2+2k,2k} = (m-1)!\ k! \
\tilde{b}_{n+w+2+2k,2k}\ .
\end{equation}
Since $\tilde{b}_{n+w+2+2k,2k} = \str(e^{-\Fs_t})_{[2k]}$ the
result can now be deduced by direct appeal to the Bismut Local
Family Index Theorem formula \cite{BGV}, Theorem(10.23), but let
us rather outline how one can deduce this from the local resolvent
symbols.

\vskip 2mm

 {\em The following computation is part of joint work with
Don Zagier \cite{SZ}}.

\vskip 2mm

For brevity we will consider the case where $\Ee$ is trivial with
zero curvature; the general case follows easily from this one.
Then from the local formula \eqref{e:smalltbsuperconnection} it is
sufficient to compute $\tilde{b}_{n+w+2+2k,2k}$ for the local
operator
\begin{equation}\label{locaH}
H =  -\sum_i \left(\dd_i - \frac{1}{4}\,r_{i}\,\xx_i \right)^2
\end{equation}
since the vertical skew-adjoint matrix of 2-forms
$(R^{M/B}\partial_i,\partial_j)$ can be written with respect to a
particular  vertical orthonormal basis as the direct sum of
$2\times 2$ blocks of 2-forms $\begin{bmatrix}
  0 & -r_{j} \\
  r_{j} & 0 \\
\end{bmatrix}$ along the fibres.
By definition one then has
\begin{equation}\label{e:heatthm}
\tilde{b}_{n+w+2+2k,2k} = \frac{1}{(2\pi)^n}\int_{M/B}
\int_{\Rf^n} \frac{i}{2\pi}\int_{C_0} e^{-\la} \ \qq_{2k}(\xx,\xi,
\la) \ d\la \ d\xi
\end{equation}
where $\qq_j(\xx,\xi, \la)$ are $2j$-forms which are the
$\xi$-homogeneous terms of the symbol of the resolvent operator
$(H - \la I)\ii$ and $C_0$ is a contour coming in on a ray with
argument in $(0,\pi/2)$, encircling the origin, and leaving on a
ray with argument in $(-\pi/2,0)$. From \eqref{locaH}, the product
formula for symbols implies that the $\qq_j$ are determined by the
following recurrence relation: let $\Delta =
\sum_{k=1}^n\dd^2_{\xx_k},$ and set $\qq_{-1} = 0$, $\qq_0 =
(|\xi|^2 - \la I)\ii$, then
\begin{equation}\label{e:recursion}
\qq_{j+1}   =   \qq_0 \, (\Delta - a_2)\, \qq_{j-1}  -  \qq_0\,
a_1\, r_j
 \hskip 10mm  (j \geq 0) \ .
\end{equation}
where with $\xi= (\xi_1,\ldots,\xi_{2n})$ $$a_1(\xx,\xi) =  i
\sum_{j=1}^{n/2}\left(\frac{r_j}{2}\right)\left(\xx_{2j-1}\,
\xi_{2j} - \xx_{2j} \, \xi_{2j-1} \right),\ \ \ \ a_2(\xx,\xi) = -
\ \frac{1}{4}
\sum_{j=1}^{n/2}\left(\frac{r_j}{2}\right)^2\left(\xx_{2j-1}^2 +
\xx_{2j}^2\right).$$ \noi It is convenient to write $\qq_j$
explicitly as a polynomial in $T = (|\xi|^2- \la I)\ii$
\begin{equation}\label{e:rjmunuT}
\qq_j(T) = \sum_{\mu+2\nu=j}\,\qq_{\mu,\nu}\,T^{\mu+\nu+1} \,,
\end{equation}
where the polynomials $\qq_{\mu,\nu}$ are now given recursively by
\begin{equation}\label{e:rjmunu}
\qq_{\mu,\nu}\=\
\left\{%
\begin{array}{ll}
    \qquad\qquad\qquad 1 &\text{if $\mu=\nu=0$,}\\
    -a_1\,\qq_{\mu-1,\nu}
  \ +(\D-a_2)\,\qq_{\mu,\nu-1} &\text{otherwise,}
\end{array}%
\right.
\end{equation}
with the convention that $\qq_{*,*}$ is to be interpreted as 0 if
either index is negative.  Then it is shown in \cite{SZ} that
\eqref{e:rjmunu} implies
\begin{equation}\label{e:Rexpansion-noXY}
\sum_{\mu,\,\nu\ge0}\frac{\qq_{\mu,\nu}(\xx,\xi)}{(\mu+\nu)!} =
\biggl(\prod_{i=1}^n\frac1{\cosh \hat{r}_i}\biggr)^{1/2}
\,\exp\biggl(|\xi|^2\, - \,\sum_{i=1}^n\,\frac{\tanh
\hat{r}_i}{\hat{r}_i}\,\bigl(\xi_i+\frac{1}{2} \,\hat{r}_i
\xx_i\bigr)^2\biggr)\ ,
\end{equation}
where $\hat{r}_{2j-1} =   -\hat{r}_{2j} = i\,r_j$. It follows that
$$\sum_{j=0}^{n/2}\int_{\Rf^n} \frac{i}{2\pi}\int_{C_0} e^{-\la} \
\qq_{2j}(\xx,\xi, \la) \ d\la \ d\xi$$
\begin{eqnarray*}\label{e:heatthmproof}
& = & \int_{\Rf^n} \sum_{\mu,\nu\geq 0}
\frac{\qq_{\mu,\nu}(\xx,\xi, \la)}{(\mu + \nu
)!}\,\underbrace{\frac{i}{2\pi}\int_{C_0} e^{-\la} \
 (|\xi|^2 - \la)\ii \ d\la}_{e^{-|\xi|^2}} \
d\xi \\[2mm]
& \stackrel{\eqref{e:Rexpansion-noXY}}{=} & 2\pi^{n/2} \,
\widehat{A}(M/B) \, .
\end{eqnarray*}
This holds for all $\xx$ and in particular at $\xx=0$. Hence
\begin{equation}\label{e:bismutlocalfamilyindex}
b_{n+w+2+2k,2k} = k!  \left((2\pi )^{-\frac{n}{2}}\int_{M/B}
\widehat{A}(M/B)
    \right)_{[2k]}   \ .
\end{equation}
Equations \eqref{e:strrescaled6}, \eqref{e:singcoeffs3},
\eqref{e:bismutlocalfamilyindex} thus combine to prove
\eqref{e:smalltlimitresstrace}.
\end{proof}

\begin{cor}
With $\Fs_t$ the Bismut superconnection curvature, the resolvent
trace differential form $\log \cc_{\, \z}(\Af_t)$ has a limit in
$\Aa(B)$ as $t\to 0$ given by the formula
\begin{equation}\label{e:smalltlimitlogzetachern}
    \lim_{t\too 0} \log \cc_{\, \z}(\Af_t)
 =  \sum_{k=0}^{[\dim B/2]}(-1)^k (k-1)!
    \, \left((2\pi )^{-\frac{n}{2}}\int_{M/B} \widehat{A}(M/B) \ch^{'}(\Ee)
    \right)_{[2k]} \ .
\end{equation}
\end{cor}
\begin{proof} We have
\begin{eqnarray}
\str(\partial^{m-1}_{\la}((\Is + \Fs_t) - \la \Is)\ii\,)  & \sim &
\sum_{k=0}^{[\dim B/2]}(-1)^{-1-k}\b_{n+w+2+2k,2k}(1-\la)^{-1-k-m}
\nonumber \\  & & + \sum_{k=0}^{[\dim B/2]}\sum_{j\geq 1}
\b_{j,2k}(1-\la)^{-\frac{j}{2}-1-k-m } \ t^{\frac{j}{2}} \ ,
\label{e:strrescaledbismutcurvature}
\end{eqnarray}
while by the previous proof \eqref{e:smalltlimitresstrace} becomes
\begin{equation}\label{e:smalltlimitresstrace2}
    \lim_{t\too 0} \str(\partial^{m-1}_{\la}((\Is + \Fs_t) - \la \Is)\ii\,)
\end{equation}
$$ = (m-1)! \sum_{k=0}^{[\dim B/2]}(-1)^{-1-k}(1-\la)^{-1-k-m}
    \, k! \left((2\pi )^{-\frac{n}{2}}\int_{M/B} \widehat{A}(M/B) \ch^{'}(\Ee)
    \right)_{[2k]} \ ,$$
and \eqref{e:singstrdet1} becomes
$$\str\left(\log(\Is + \Fs_t)\,(\Is + \Fs_t)\si \right)|^{{\rm mer}}
$$
$$
=  \sum_{k=0}^{[\dim B/2]}\, (-1)^k k! \left((2\pi
)^{-\frac{n}{2}}\int_{M/B} \widehat{A}(M/B) \ch^{'}(\Ee)
    \right)_{[2k]} \ \dd_s F_{1+k}(s-1)|^{{\rm mer}} $$
\begin{equation}\label{e:singstrdet2}
+ \ \dd_s h_{2k+3+n+w,t}(s) \ ,
\end{equation}
where $h_{N,t}(s)$ is the remainder term $h_{N}(s)$ for the
$t$-rescaled Bismut superconnection with $N = 2k+3+n+w$. The
coefficients in the above formulas are precisely related as
explained in the previous proof. The powers in the expansion
\eqref{e:strrescaledbismutcurvature} have the crucial consequence
for the globally determined remainder term in
\eqref{e:singstrdet2} that
\begin{equation}\label{e:scalestozero}
\dd_s h_{2k+3+n+w,t}(s) = O(t^{1/2})
\end{equation}
as $t\to 0+$. It remains then to compute the sum in
\eqref{e:singstrdet2}. We have for $-1-\re(\a) < \re(s) < 0$
\begin{eqnarray}
\dd_s F_{1+\a}(s-1) & = & \dd_s \, \frac{i}{2\pi}\int_C \mu^{-s}
(1-\mu)^{-\a-1} \ d\mu \nonumber\\
& = & -\frac{1}{
\a}\, \frac{i}{2\pi}\int_C \log\mu \,\mu^{-s}
\dd_{\mu}(1-\mu)^{-\a} \ d\mu \nonumber\\
& = & \frac{1}{\a}\, \frac{i}{2\pi}\int_C \mu^{-s-1} (1-\mu)^{-\a}
\ d\mu - \frac{s}{\a}\, \frac{i}{2\pi}\int_C \log\mu
\,\mu^{-s-1}(1-\mu)^{-\a} \ d\mu \nonumber\\
& = & \frac{1}{\a} F_{\a}(s) - \frac{s}{\a} \,\dd_s\, F_{\a}(s) \
. \label{e:Falphas}
\end{eqnarray}
Since $\G(z)|^{{\rm mer}}$ is holomorphic for $z\neq 0, -1,
-2,\ldots$, then from \eqref{e:Ftsidentities} we have that
$F_{\a}(s)|^{{\rm mer}}$ is holomorphic near $s=0$ and hence that
\begin{equation*}
\left(\frac{s}{\a} \,\dd_s\, F_{\a}(s)\right)|^{{\rm mer}}_{s=0} =
0 \ .
\end{equation*}
The meromorphic extension of the Gamma function is obtained by the
identity $\G(s+1)=s\G(s)$ ,and hence from \eqref{e:Ftsidentities}
we find that
\begin{equation*}
F_{k}(s)|^{{\rm mer}}_{s=0} = 1 \ .
\end{equation*}
Hence
\begin{equation*}
\dd_s F_{1+k}(s-1)|^{{\rm mer}}_{s=0} = \frac{1}{k} \ .
\end{equation*}
From \eqref{e:singstrdet2} and \eqref{e:scalestozero} we obtain
the asserted identity
\end{proof}

\vskip 5mm

Using \propref{p:resformatinfinity}, one has for the Bismut
superconnection that the limit
$$\lim_{\e\to 0}\int_{\e}^{1/\e} \zeta_{\pi}((\Is +
\Fs_{\e})\ii\dot{\Af}_{\e}, \Is + \Fs_{\e},0)|^{{\rm mer}} \,
d\e$$ exists uniformly in all $C^l$ norms on compact subsets of
$B$, and hence that
$$\tau_{0,\o}(\Af) := \lim_{\e\to 0}\tau_{\e,\e\ii}(\Af)$$
exists. With \propref{p:zetacherntransgression2}, this completes
the proof of the local family index formula for the zeta-Chern
class:
\begin{eqnarray*}
\log \cc(\Ker (\Ds),\nabla^0)  & = &  \sum_{k=0}^{[\dim
B/2]}(-1)^k (k-1)!
    \, \left((2\pi )^{-\frac{n}{2}}\int_{M/B} \widehat{A}(M/B) \ch^{'}(\Ee)
    \right)_{[2k]} \\
& & +  \   d \ \lim_{\e\to 0}\int_{\e}^{1/\e} \zeta_{\pi}((\Is +
\Fs_{\e})\ii\dot{\Af}_{\e}, \Is + \Fs_{\e},0)|^{{\rm mer}} \, d\e
\ ,
\end{eqnarray*}
or, exponentiating,
\begin{equation*}
\cc(\Ker (\Ds),\nabla^0) =  \prod_{k=0}^{[\dim B/2]}e^{(-1)^k
(k-1)!
    \, \left((2\pi )^{-\frac{n}{2}}\int_{M/B} \widehat{A}(M/B) \ch^{'}(\Ee)
    \right)_{[2k]}} +    d\omega_{0,\o} \ .
\end{equation*}

\vskip 8mm

\noi {\small King's College London \\ simon.scott@kcl.ac.uk}


\begin{thebibliography}{99}

\bibitem[BF]{BF} Bismut, J-M and Freed, D.S.: 1986,
`The analysis of elliptic families I', {\it Comm. Math. Phys. {\bf
106}}, 159--176.

\bibitem[BGV]{BGV}  Berline, N., E. Getzler, and M. Vergne: {\it Heat
Kernels and Dirac Operators}.  Grundlehren der Mathematischen
Wissenschaften {\bf 298}, Springer-Verlag, Berlin, 1992.

\bibitem[B]{B} Bismut, J-M.: 1986, `The Atiyah-singer index theorem for
families of Dirac operators: Two heat equation proofs', {\it
Invent. Math {\bf 83}}, 91--151.

\bibitem[Gr1]{Gr1} G. Grubb, `A resolvent approach to
traces and zeta Laurent expansions', AMS Contemp. Math. Proc.,
vol. 366, 67--93 (2005). See also  arXiv: math.AP/0311081.

\bibitem[GH]{GH} Grubb, G., Hansen, L.: 2002, `Complex powers of
resolvents of pseudodifferential operators', {\it Comm. Part.
Diff. Eq. {\bf 27}}, 2333-2361.

\bibitem[GS1]{GS1} Grubb, G., Seeley, R.: 1995, `Weakly parametric
pseudodifferential operators and Atiyah-Patodi-Singer boundary
problems', {\it Invent. Math. {\bf 121}}, 481-529.

\bibitem[GS2]{GS2} Grubb, G., Seeley. R.: 1996,
`Zeta and eta functions for Atiyah-Patodi-Singer operators', {\it
J. Geom. Anal. {\bf 6}}, 31--77.

\bibitem[Q]{Q} Quillen, D.G.: 1985, `Superconnections and the Chern character',
{\it Topology {\bf 24}}, 89--95.

\bibitem[S]{S} Seeley, R. T.: 1967, `Complex powers of an elliptic
operator', AMS {\it Proc. Symp. Pure Math. X}. AMS Providence,
288--307.

\bibitem[SZ]{SZ} Scott, S., Zagier. D: `A symbol proof of the local Atiyah-Singer index theorem', in preparation.

\bibitem[Sh]{Sh} Shubin, M.A.: {\it Pseudodifferential Operators and Spectral Theory', 2nd
Edition}, Springer, 2001.



\end{thebibliography}
\end{document}